# LOCAL LINEAR QUANTILE ESTIMATION FOR NONSTATIONARY TIME SERIES

By Zhou Zhou and Wei Biao Wu

*University of Chicago*

We consider estimation of quantile curves for a general class of nonstationary processes. Consistency and central limit results are obtained for local linear quantile estimates under a mild short-range dependence condition. Our results are applied to environmental data sets. In particular, our results can be used to address the problem of whether climate variability has changed, an important problem raised by IPCC (Intergovernmental Panel on Climate Change) in 2001.

**1. Introduction.** Nonstationary time series analysis has attracted considerable attention recently; see, for example, Priestley (1988), Dahlhaus (1997), Mallat, Papanicolaou and Zhang (1998), Nason, von Sachs and Kroisandt (2000), Giurcanu and Spokoiny (2004) and Ombao, von Sachs and Guo (2005), among others. Many of the previous results concern estimation of second-order characteristics such as covariance functions and time-varying spectra.

The primary goal of the paper is to estimate quantile curves of a special class of nonstationary processes that can be loosely called locally stationary processes. Conceptually, local stationarity means that the physical mechanism generating such processes changes smoothly in time [Mallat, Papanicolaou and Zhang (1998)]. For the purpose of estimating quantile curves of locally stationary processes, the classical framework which is based on second-order characteristics does not seem to be quite suitable. In particular, to estimate quantile functions of heavy-tailed processes, one cannot directly apply results which are based on covariance functions.

Here we shall adopt the following formulation. Let the observed sequence $\{X_{i,n}\}_{i=1}^{n}$ be generated from the model

$$X_{i,n} = G(i/n, \Upsilon_i), \tag{1}$$

where $\Upsilon_i = (\ldots, \varepsilon_{i-1}, \varepsilon_i)$ and $\varepsilon_j$, $j \in \mathbb{Z}$ are independent and identically distributed (i.i.d.) random variables, and $G : [0,1] \times \mathbb{R}^\infty \mapsto \mathbb{R}$ is a measurable function such that $\zeta_i(t) := G(t, \Upsilon_i)$ is a properly defined random variable for all $t \in [0,1]$. Let $F(t,x) = \mathbb{P}(\zeta_i(t) \leq x)$, $x \in \mathbb{R}$ be the cumulative distribution function (c.d.f.) of $\zeta_i(t)$, $t \in [0,1]$; let $Q_\alpha(t)$ be the $\alpha$th quantile function of $\zeta_i(t)$, $\alpha \in (0,1)$, namely $Q_\alpha(t) = \inf_x \{F(t,x) \geq \alpha\}$. We shall estimate $Q_\alpha(t)$, the $\alpha$th quantile curve, from $\{X_{i,n}\}_{i=1}^n$ by the local linear approach. In the sequel for notational convenience we shall write $X_{i,n}$ as $X_i$.

Process (1) is nonstationary. It covers a wide range of nonstationary processes and it naturally extends many existing stationary time series models into the nonstationary setting; see examples in Section 4. Wiener (1958) claimed that, for any stationary ergodic process $(X_i)$, there exists a measurable function $G$ and i.i.d. random variables $\varepsilon_i$, $i \in \mathbb{Z}$, such that $X_i = G(\Upsilon_i)$ holds for all $i$; see also Tong (1990), page 204. The latter relation can be interpreted as a physical system with $\Upsilon_i$ and $X_i$ being the input and the output, respectively, and $G$ being the transform or filter that represents the underlying physical mechanism. See Wu (2005) for further references. In (1), by allowing the data generating mechanism $G$ depending on the time index $t$ in such a way that $G(t, \Upsilon_i)$ changes smoothly with respect to $t$, one has local stationarity in the sense that the subsequence $\{X_i, \ldots, X_{i+j-1}\}$ is approximately stationary if its length $j$ is sufficiently small compared to $n$.

Estimation of quantile curves is an important problem. The report *Climate Change 2001* written by the Intergovernmental Panel on Climate Change (IPCC) addresses the serious problem of whether the climate variability or climate extremes have changed. Prominent changes in climate variability can have a devastating impact upon the environment. In this paper, we shall interpret climate extremes as upper and lower quantiles such as $Q_{0.95}(t)$ or $Q_{0.05}(t)$. By estimating such quantile curves, we obtain more detailed distributional information about the underlying process than the mean and the variability based on covariances; see Yu, Lu and Stander (2003) and Koenker (2005) for more discussions.

There is a large literature on state domain nonparametric conditional quantile estimation where an i.i.d. or stationary series $(Y_i, X_i)_{i=1}^n$ is observed and one is interested in estimating the conditional quantiles of $Y_i$ given $X_i$; see, for instance, Bhattacharya and Gangopadhyay (1990), Chaudhuri (1991), Chapter 6 of Fan and Gijbels (1996) and Yu and Jones (1998) for i.i.d. data and Abberger (1997), Abberger and Heiler (2002), Cai (2002) and Cai and Xu (2005) for stationary strong mixing time series. For other contributions see Cheng and Parzen (1997) and Csörgő and Mielniczuk (1996). On the other hand, there are much fewer results on time domain nonparametric quantile estimation, and people are mostly interested in the subordinated Gaussian process $X_i = H(i/n, Z_i)$, $i = 1, 2, \ldots, n$, where $(Z_i)$ is a stationary

Gaussian process and $H$ is a measurable function; see Ghosh, Beran and Innes (1997) and Ghosh and Draghicescu (2002a, 2002b), among others.

The rest of the paper is structured as follows. Section 2 introduces the local linear quantile estimates, dependence measures and regularity conditions. Central limit theorems and Bahadur representations are stated in Section 3 and proved in Section 6. Section 4 presents examples of nonstationary linear processes and nonstationary nonlinear time series models. Applications to environmental data are given in Section 5.

**2. Preliminaries.** We now introduce some notation. For a vector $\mathbf{v} = (v_1, v_2, \ldots, v_p) \in \mathbb{R}^p$, let $|\mathbf{v}| = (\sum_{i=1}^p v_i^2)^{1/2}$. For a $p \times p$ matrix $A$, define $|A| = \sup\{|A\mathbf{v}| : |\mathbf{v}| = 1\}$. For a random vector $\mathbf{V}$, write $\mathbf{V} \in \mathcal{L}^q$ $(q > 0)$ if $\|\mathbf{V}\|_q := [\mathbb{E}(|\mathbf{V}|^q)]^{1/q} < \infty$ and $\|\mathbf{V}\| = \|\mathbf{V}\|_2$. Denote by $\Rightarrow$ the weak convergence. For an interval $\mathcal{I} \subset \mathbb{R}$, denote by $\mathcal{C}^i \mathcal{I}$, $i \in \mathbb{N}$, the collection of functions that have $i$th order continuous derivatives on $\mathcal{I}$, and, for $\mathcal{D} \subset \mathbb{R}^d$, let $\mathcal{C}\mathcal{D}$ be the collection of real-valued functions that are continuous on $\mathcal{D}$. A function $f : \mathbb{R}^d \to \mathbb{R}$ is said Lipschitz continuous on $\mathcal{D} \subset \mathbb{R}^d$ if there exists a finite constant $C$, such that $|f(x_1) - f(x_2)| \leq C|x_1 - x_2|$ for all $x_1, x_2 \in \mathcal{D}$. For $x \in \mathbb{R}$, define $x^+ = \max(x, 0)$. The symbol $C$ denotes a finite generic constant which may vary from line to line.

2.1. *Local linear quantile estimator.* Observe that $Q_\alpha(t)$, the $\alpha$th quantile function of $\zeta_i(t)$, satisfies

$$(2) \qquad Q_\alpha(t) = \arg\min_{\beta_0} \mathbb{E}\{\rho_\alpha(\zeta_i(t) - \beta_0) - \rho_\alpha(\zeta_i(t))\},$$

where $\rho_\alpha(x) = \alpha x^+ + (1-\alpha)(-x)^+$ is the check function [Koenker (2005)]. As $Q_\alpha(t_1) \approx Q_\alpha(t) + (t_1 - t)Q'_\alpha(t)$ for $t_1$ close to $t$, it is natural to estimate $Q_\alpha(t)$ and $Q'_\alpha(t)$ via

$$(3)\ (\hat{Q}_{\alpha,b_n}(t), \hat{Q}'_{\alpha,b_n}(t)) = \arg\min_{(\beta_0, \beta_1)} \sum_{i=1}^n \rho_\alpha(X_i - \beta_0 - \beta_1(i/n - t))K_{b_n}(i/n - t),$$

where $K$ is a kernel function, $K_{b_n}(\cdot) = K(\cdot/b_n)$ and $b_n > 0$ is the bandwidth. We shall omit the subscript $b_n$ in $\hat{Q}$ and $\hat{Q}'$ hereafter if no confusion will be caused.

Note that (3) defines a local linear quantile estimate, while the analogous version of (2) gives a local constant estimate; see (10). Draghicescu, Guillas and Wu (2008) considered local constant quantile estimation. In the context of local polynomial regression based on mean squared errors (MSE), Fan and Gijbels (1996) argued that the local linear estimate is generally better than the local constant estimate since the latter suffers from the notorious boundary problem. Proposition 1 below asserts that, for quantile curve estimation, the local linear approach can also alleviate the boundary problem.

2.2. *Dependence measures and regularity conditions.* To study $\hat{Q}$ in (3), it is necessary to introduce appropriate dependence measures. Following Wu (2005), we shall adopt *physical dependence measures* for a family of stochastic processes. All our results will be expressed in terms of these dependence measures.

DEFINITION 1. Let $\{H(t, x, \Upsilon_i)\}_{i \in \mathbb{Z}}$, $(t, x) \in [0, 1] \times \mathbb{R}$, be a family of stochastic processes with $H(t, x, \Upsilon_i) \in \mathcal{L}^2$. Let $\varepsilon_0', \varepsilon_j$, $j \in \mathbb{Z}$, be i.i.d. For $k \geq 0$, let $\Upsilon_k^* = (\Upsilon_{-1}, \varepsilon_0', \varepsilon_1, \ldots, \varepsilon_k)$ and define the physical dependence measure

$$(4) \qquad \delta_H(k) = \sup_{(t,x) \in [0,1] \times \mathbb{R}} \|H(t, x, \Upsilon_k) - H(t, x, \Upsilon_k^*)\|.$$

Here we recall $\|\cdot\| = [\mathbb{E}(|\cdot|^2)]^{1/2}$. We say that the system $\{H(\cdot, \cdot, \Upsilon_i)\}_{i \in \mathbb{Z}}$ is stable if

$$(5) \qquad \sum_{k=0}^{\infty} \delta_H(k) < \infty.$$

Observe that $\delta_H(k)$ measures the dependence of $H(t, x, \Upsilon_k)$ on the input $\varepsilon_0$ over $(t, x) \in [0, 1] \times \mathbb{R}$. Condition (5) indicates that the cumulative impact of $\varepsilon_0$ on future values is bounded, thus suggesting short-range dependence. At a fixed $(t, x) \in [0, 1] \times \mathbb{R}$, if $H(t, x, \Upsilon_i)$ is a linear process, (5) reduces to the stability condition introduced in Box, Jenkins and Reinsel (1994). All the results of this paper will be derived under this short-range dependence condition. If the stability condition (5) is violated, namely if the sum in equation (5) is infinite, then we have a long-memory process and the asymptotic theory developed in the paper does not work. Note that $\delta_H(k)$ is closely related to the data generating mechanism and hence is easy to work with; see examples in Section 4 and also Wu (2005).

Hereafter we shall assume that $K(\cdot) \in \mathcal{K}$, where $\mathcal{K}$ is the collection of density functions $K$ such that $K$ is symmetric with support $[-1, 1]$ and $K \in \mathcal{C}^1[-1, 1]$. A popular choice is the Epanechnikov kernel $K(u) = 3\max(0, 1 - u^2)/4$. For $K \in \mathcal{K}$ let $\phi_K = \int_{-1}^{1} K^2(u)\,du$. Define $\mu_{j,K} = 2\int_0^1 u^j K(u)\,du$ if $j \in \mathbb{N}$ is even, and $\mu_{j,K} = \int_0^1 u^j K(u)\,du$ if $j \in \mathbb{N}$ is odd. We shall abbreviate $\phi_K$ (resp. $\mu_{j,K}$) as $\phi$ (resp. $\mu_j$) if no confusion will be caused.

Note that, for fixed $t$, $\{\zeta_i(t)\}_{i \in \mathbb{Z}}$ is a stationary process. We shall impose the following regularity conditions:

(A1) For every $\alpha \in (0, 1)$, $Q_\alpha(\cdot) \in \mathcal{C}^2[0, 1]$.

(A2) $f(t, x) \in \mathcal{C}([0, 1] \times \mathbb{R})$, where $f(t, \cdot)$ is the density function of $\zeta_i(t)$.

(A3) (Stochastic Lipschitz Continuity) There exists $C, q > 0$, such that $\|\zeta_i(t_1) - \zeta_i(t_2)\|_q \leq C|t_1 - t_2|$ holds for all $t_1, t_2 \in [0, 1]$.

Condition (A3) indicates that the underlying data generating mechanism $G(t, \cdot)$ changes smoothly in time, thus suggesting local stationarity. Conditions (A1)–(A3) are needed in the sense that, without them, it might be impossible to estimate quantile curves of the series from its stochastic variation [see the discussions in Chapter 6 in Fan and Yao (2003)]. We will not consider the possibility of jumps in quantile curves in this paper.

## 3. Main results.

3.1. *Asymptotic normality.* Let $\psi_\alpha(x) = \alpha - I\{x \leq 0\}$ be the left derivative of $\rho_\alpha(x)$. Theorem 1 asserts that, for $t \in (0,1)$, $nb_n f(t, Q_\alpha(t))[\hat{Q}_\alpha(t) - Q_\alpha(t)]$ can be approximated by the linear form

$$T_{\alpha,n}(t) := \sum_{i=1}^{n} \psi_\alpha(X_i - Q_\alpha(t) - Q'_\alpha(t)(i/n - t))K_{b_n}(i/n - t)$$

with an $o_\mathbb{P}((nb_n)^{1/2})$ error. Due to the linearity it is easier to deal with $T_{\alpha,n}(t) - \mathbb{E}[T_{\alpha,n}(t)]$ which is asymptotically normal under proper conditions. Recall $\phi = \int_{-1}^{1} K^2(u)\,du$ and $\zeta_i(t) = G(t, \Upsilon_i)$.

THEOREM 1. *Assume $b_n \to 0$ and $nb_n \to \infty$. Let $t \in (0,1)$, $J(t, x, \Upsilon_i) = I\{G(t, \Upsilon_i) \leq x\}$ and $F(t, x, \Upsilon_i) = \mathbb{P}(G(t, \Upsilon_{i+1}) \leq x | \Upsilon_i)$. Further assume $f(t, Q_\alpha(t)) > 0$, (A1) and (A2), and that the family $\{F(\cdot, \cdot, \Upsilon_i)\}_{i \in \mathbb{Z}}$ is stable. Then*

(6) $\quad (nb_n + O(1))f(t, Q_\alpha(t))[\hat{Q}_\alpha(t) - Q_\alpha(t)] - T_{\alpha,n}(t) = o_\mathbb{P}((nb_n)^{1/2}).$

*Additionally, assume* (A3), $nb_n^5 = O(1)$ *and*

(7) $$\sigma^2(t) := \left\| \sum_{i=0}^{\infty} \mathcal{P}_0 J(t, Q_\alpha(t), \Upsilon_i) \right\|^2 > 0,$$

*where $\mathcal{P}_0(\cdot) = \mathbb{E}(\cdot | \Upsilon_0) - \mathbb{E}(\cdot | \Upsilon_{-1})$. Then the Central Limit Theorem (CLT) holds:*

(8) $\quad (nb_n)^{1/2}[\hat{Q}_\alpha(t) - Q_\alpha(t) - b_n^2 \mu_2 Q''_\alpha(t)/2] \Rightarrow \mathrm{N}(0, \phi\sigma^2(t)/f^2(t, Q_\alpha(t))).$

Since $\{F(\cdot, \cdot, \Upsilon_i)\}_{i \in \mathbb{Z}}$ is stable, $\sigma^2(t)$ exists and is finite; see (37). We call $\sigma^2(t)$ the long-run variance which is due to the dependence of the series. If $J(t, Q_\alpha(t), \Upsilon_i), j \in \mathbb{Z}$, are independent, then $\sigma^2(t) = \mathrm{var}(J(t, Q_\alpha(t), \Upsilon_i)) = \alpha(1-\alpha)$. Theorem 1 entails the following corollary on interquartile range (IQR), a simple measure for the spread of distributions.

COROLLARY 1. *Let $IQR(t) = Q_{0.75}(t) - Q_{0.25}(t)$ and $\widehat{IQR}(t) = \hat{Q}_{0.75}(t) - \hat{Q}_{0.25}(t)$. Let conditions in Theorem 1 be satisfied with $\alpha = 0.75$ and $\alpha = 0.25$. Let $J_{IQR}(t, \Upsilon_i) = J(t, Q_{0.75}(t), \Upsilon_i)/f(t, Q_{0.75}(t)) - J(t, Q_{0.25}(t), \Upsilon_i)/f(t, Q_{0.25}(t))$ and assume that $\sigma_{IQR}(t) := \|\sum_{i=0}^{\infty} \mathcal{P}_0 J_{IQR}(t, \Upsilon_i)\| > 0$. Then*

$$(9) \quad (nb_n)^{1/2}[\widehat{IQR}(t) - IQR(t) - b_n^2 \mu_2 IQR''(t)/2] \Rightarrow N(0, \phi \sigma_{IQR}^2(t)).$$

An important issue in applying Theorem 1 is the choice of the bandwidth $b_n$. A data-driven procedure is given in Section 3.1.1.

We now study the boundary behavior of the local linear estimates. In the context of quantile estimation for stationary processes, Cai and Xu (2005) argued that the local polynomial estimates have a better performance than the local constant estimates. Proposition 1 concerns left boundary points. Analogous results hold for right boundary points. It shows that the local linear quantile estimator alleviates the boundary problem in that the asymptotic MSEs induced by (8) and (12) are of the same order. For local constant quantile estimate, however, (11) implies that its rate of convergence is slower at the boundary points. Let $Q'_\alpha(0+), Q''_\alpha(0+), \ldots$, be right derivatives of $Q_\alpha$ at 0.

PROPOSITION 1 (Boundary behavior). *Define the local constant estimate*

$$(10) \quad \bar{Q}_{\alpha, b_n}(t) = \arg\min_{\beta} \sum_{i=1}^{n} \rho_\alpha(X_i - \beta) K_{b_n}(i/n - t).$$

*Under assumptions of Theorem 1 with $t = 0$, we have:* (i)

$$(11) \quad \begin{aligned} &(nb_n/2)^{1/2} f(0, Q_\alpha(0))\{\bar{Q}_\alpha(0) - Q_\alpha(0) - 2b_n[Q'_\alpha(0+)\mu_1 + o(1)]\} \\ &\Rightarrow N(0, \sigma^2(0)\phi) \end{aligned}$$

*and* (ii) *for the local linear estimate $\hat{Q}_\alpha(0)$,*

$$(12) \quad \begin{aligned} &(nb_n)^{1/2} f(0, Q_\alpha(0))[\hat{Q}_\alpha(0) - Q_\alpha(0) - b_n^2 Q''_\alpha(0+) B_K/2] \\ &\Rightarrow N(0, \sigma^2(0) V_K), \end{aligned}$$

*where $B_K = (\mu_2^2 - 4\mu_1\mu_3)/(\mu_2 - 4\mu_1^2)$ and $V_K = 4\int_0^1 (\mu_2 - 2\mu_1 u)^2 K(u)^2 \, du / (\mu_2 - 4\mu_1^2)^2$.*

3.1.1. *Bandwidth selection.* For the estimate $\hat{Q}_\alpha(\cdot) = \hat{Q}_{\alpha, b_n}(\cdot)$, based on (8) of Theorem 1, we define the weighted asymptotic mean integrated squared error (AMISE) as

$$AMISE(\hat{Q}_{\alpha, b_n}) = \int_0^1 \{[b_n^2 \mu_2 Q''_\alpha(t)/2]^2 + [\phi \sigma^2(t)/(nb_n)]\} f^2(t, Q_\alpha(t)) \, dt.$$

Here the weight $f^2(t, Q_\alpha(t))$ indicates the abundance of data for estimating $Q_\alpha(t)$. With more data in the neighborhood of $(t, Q_\alpha(t))$, we can have a more accurate estimate and therefore a larger weight is assigned. Observe that the optimal AMISE bandwidth

$$(13) \qquad b_n(\alpha) = \left\{ \frac{\phi \int_0^1 \sigma^2(t)\,dt}{n\mu_2^2 \int_0^1 [Q''_\alpha(t)]^2 f^2(t, Q_\alpha(t))\,dt} \right\}^{1/5}$$

is of order $n^{-1/5}$. For independent data, several authors have considered the bandwidth selection problem for quantile smoothing; see, for example, Fan and Gijbels (1996) and Yu and Jones (1998) (YJ hereafter). Denote by $b_n^{\text{ind}}(\alpha)$ the optimal bandwidth obtained under independence. By YJ and (13),

$$(14) \qquad \frac{b_n(\alpha)}{b_n^{\text{ind}}(\alpha)} = \left[ \frac{\int_0^1 \sigma^2(t)\,dt}{\alpha(1-\alpha)} \right]^{1/5} := \rho^*(\alpha).$$

Note that $\alpha(1-\alpha) = \text{var}(J(t, Q_\alpha(t), \Upsilon_i))$ and $\rho^*(\alpha)$ is called the variance correction factor which accounts for the dependence. YJ proposed a simple and easy-to-use rule of thumb selector for $b_n^{\text{ind}}(\alpha)$, which we denote by $b_{n,YJ}(\alpha)$. Our idea is to use $b_n(\alpha)$ by correcting $b_{n,YJ}(\alpha)$ by a factor of $\hat{\rho}^*(\alpha)$, an estimate of $\rho^*(\alpha)$. More precisely, let

$$\hat{\rho}^*(\alpha) = \frac{(\tilde{\sigma}^2)^{1/5}}{(\alpha(1-\alpha))^{1/5}},$$

where

$$\tilde{\sigma}^2 = \frac{\tilde{m}}{n - \tilde{m} + 1} \sum_{j=1}^{n-\tilde{m}+1} \left( \frac{1}{\tilde{m}} \sum_{i=j}^{j+\tilde{m}-1} \varsigma_{i,\alpha} - \bar{\varsigma}_{n,\alpha} \right)^2,$$

$\varsigma_{i,\alpha} = \psi_\alpha(X_i - \hat{Q}_{\alpha, b_{n,YJ}(\alpha)}(t))$, $\bar{\varsigma}_{n,\alpha} = \sum_{i=1}^n \varsigma_{i,\alpha}/n$ and $\tilde{m} = \lfloor n^{1/3} \rfloor$. Here $\tilde{\sigma}^2$ is a block estimate of long-run variance. See also discussions in Section 3.4. Following similar arguments as those of Theorem 5, it can be shown that $\tilde{\sigma}^2$ converges weakly to $\int_0^1 \sigma^2(t)\,dt$. We suggest

$$(15) \qquad b_n^*(\alpha) = b_{n,YJ}(\alpha) \times \hat{\rho}^*(\alpha).$$

As pointed out by a referee, there are time blocks where the distribution of the time series (or a few particular quantiles of interest) changes more quickly than others. This suggests using a bandwidth that could change with time. In general, variable bandwidth selection can be done by performing (15) locally. More precisely, similar to (14), the optimal asymptotic MSE local bandwidth $b_n(\alpha, t)$ satisfies

$$\frac{b_n(\alpha, t)}{b_n^{\text{ind}}(\alpha, t)} = \left[ \frac{\sigma^2(t)}{\alpha(1-\alpha)} \right]^{1/5} := \rho^*(\alpha, t),$$

where $b_n^{\mathrm{ind}}(\alpha, t)$ is the optimal MSE local bandwidth under independence and $\rho^*(\alpha, t)$ is called the local variance correction factor. Here $\sigma^2(t)$ can be estimated using (25) in Section 3.4.

3.2. *Bahadur representations.* Bahadur representation is an important tool for asymptotic analysis of estimators and it provides deep insight into properties of estimators by approximating them by linear forms. See, for instance, He and Shao (1996), Koenker (2005) and Wu (2007a), and references therein. In this section we shall obtain both locally and globally uniform Bahadur representations for the estimated quantile curve $\hat{Q}_\alpha(\cdot)$. Define $\hat{\theta}_{\alpha,n}(t) = (\hat{\theta}_{\alpha,n,1}(t), \hat{\theta}_{\alpha,n,2}(t))^\top := (\hat{Q}_\alpha(t) - Q_\alpha(t), b_n(\hat{Q}'_\alpha(t) - Q'_\alpha(t)))^\top$. A key component to derive Bahadur representation for $\hat{\theta}_{\alpha,n}(s)$ is to study the oscillation behavior of the associated weighted empirical process of form (16). Recall that Bahadur (1966) also analyzed oscillations of empirical processes and obtained asymptotic expansions of sample quantiles of i.i.d. data.

Let $\mathbf{z}_{i,n}(t) = (1, (i/n - t)/b_n)^\top$. For $\theta = (\theta_1, \theta_2)^\top \in \mathbb{R}^2$, define

$$S_{\alpha,n}(t, \theta) = \sum_{i=1}^n \psi_\alpha(X_i - Q_\alpha(t) - (i/n - t)Q'_\alpha(t) - \theta^\top \mathbf{z}_{i,n}(t)) \tag{16}$$

$$\times K_{b_n}(i/n - t)\mathbf{z}_{i,n}(t).$$

Write $S_{\alpha,n}(t) = S_{\alpha,n}(t, \mathbf{0})$, where $\mathbf{0} = (0, 0)^\top$. Let $\mu_\mathbf{K} = \mathrm{diag}(1, \mu_2)$ be a $2 \times 2$ diagonal matrix. We need the following regularity conditions:

(B1) Assume that, for $k = 0, 1, 2, 3$, $F_k(t, x, \Upsilon_i) := \partial^k F(t, x, \Upsilon_i)/\partial x^k$ exists and the family $\{F_k(\cdot, \cdot, \Upsilon_i)\}_{i \in \mathbb{Z}}$ is stable, where $F_0(t, x, \Upsilon_i) = F(t, x, \Upsilon_i)$.

(B2) There exists $C_0 < \infty$ such that $\sup_{(t,x) \in [0,1] \times \mathbb{R}} F_1(t, x, \Upsilon_i) < C_0$ almost surely.

(B3) $f(t, x)$ is Lipschitz continuous on $[0, 1] \times \mathbb{R}$.

THEOREM 2 (Local uniform Bahadur representation). *Let $0 < t < 1$. Assume that $f(t, Q_\alpha(t)) > 0$, $nb_n/\log^2 n \to \infty$ and $b_n \to 0$. Then under* (A1), (B1)–(B3), *we have*

$$\sup_{s \in B_n(t)} \left| f(s, Q_\alpha(s)) \mu_\mathbf{K} \hat{\theta}_{\alpha,n}(s) - \frac{S_{\alpha,n}(s)}{nb_n} \right| = O_\mathbb{P}\left( \frac{\varrho_n^{1/2} \log n}{\sqrt{nb_n}} + b_n \varrho_n + \varrho_n^2 \right), \tag{17}$$

*where $B_n(t) = [t - b_n, t + b_n] \cap [0, 1]$ and $\varrho_n = (nb_n)^{-1/2} \log n + b_n^2$.*

THEOREM 3 (Global uniform Bahadur representation). *Let $T_n = [b_n, 1 - b_n]$. Assume* (A1), (B1)–(B3), $\inf_{t \in [0,1]} f(t, Q_\alpha(t)) > 0$, $nb_n^2 \to \infty$ *and $b_n \to*

0. *Then*

$$(18) \quad \sup_{t \in T_n} \left| f(t, Q_\alpha(t))\mu_{\mathbf{K}}\hat{\theta}_{\alpha,n}(t) - \frac{S_{\alpha,n}(t)}{nb_n} \right| = O_{\mathbb{P}}\left( \frac{\pi_n^{1/2} \log n}{\sqrt{nb_n}} + b_n \pi_n + \pi_n^2 \right),$$

*where* $\pi_n = (nb_n)^{-1/2}(\log n + (b_n)^{-1/2} + (nb_n^5)^{1/2})$.

Bahadur representations in Theorems 2 and 3 are useful for studying asymptotic properties of the estimated quantile curves by providing asymptotic approximation of $\hat{\theta}_{\alpha,n}(t)$. The local Bahadur representation of Theorem 2 is needed in proving consistency and asymptotic normality for a two-step smoother given in Section 3.3. By Remark 1 in Section 6, Theorem 3 implies global uniform consistency of the estimated quantile curves. We believe that, with the uniform Bahadur representation (18), one can construct simultaneous confidence bands for $Q_\alpha(t)$ over $t \in T_n$ based on a Gaussian approximation of $\{S_{\alpha,n}(t), t \in T_n\}$; see Wu and Zhao (2007) for simultaneous inference of mean trends in time series.

3.3. *Two-stage smoothing and jackknife.* The estimated curve $\hat{Q}_\alpha(t)$ may not be smooth and thus is visually unattractive. To remedy this problem, we suggest a second-stage local linear smoothing

$$(19) \quad \check{Q}_\alpha(t) = \sum_{i=1}^n \hat{Q}_\alpha(i/n) w_n(t, i),$$

where $w_n(t, i) = K_{\bar{b}_n}(t - i/n)[\mathcal{B}_2(t) - (t - i/n)\mathcal{B}_1(t)]/[\mathcal{B}_2(t)\mathcal{B}_0(t) - \mathcal{B}_1^2(t)]$ are the local linear weights, $\mathcal{B}_j(t) = \sum_{i=1}^n (t - i/n)^j K_{\bar{b}_n}(t - i/n)$ and $\bar{b}_n$ is another bandwidth; see Fan and Gijbels (1996), page 20, for the validation of the local linear weights. Fan and Zhang (2000) applied a two-step procedure to improve the unsmooth raw estimates. Introducing another smoothing step may possibly bring extra bias and variance. However, fortunately, if the bandwidth $\bar{b}_n$ satisfies $\bar{b}_n/b_n \to 0$, then the extra bias and variance are negligible compared to those of the original quantile estimator $\hat{Q}_\alpha(t)$. The simulation study in Draghicescu, Guillas and Wu (2008) suggests that the second-stage smoothing can make the raw estimate visually attractive. However the latter paper does not provide theoretical justification of the validity of the procedure. Theorem 4 presents properties of the estimator (19). The key tool for proving the theorem is the local Bahadur representation (17).

THEOREM 4. *Assume* $\bar{b}_n/b_n \to 0$, $n\bar{b}_n \to \infty$ *and* $nb_n^5 = O(1)$. *Under conditions of Theorem 2, we have*

$$(20) \quad \sqrt{n\bar{b}_n} f(t, Q_\alpha(t))[\check{Q}_\alpha(t) - Q_\alpha(t)] - T_{\alpha,n}(t)/\sqrt{n\bar{b}_n} = o_{\mathbb{P}}(1).$$

*Hence*

$$\sqrt{nb_n}[\check{Q}_\alpha(t) - Q_\alpha(t) - b_n^2\mu_2 Q''_\alpha(t)/2] \Rightarrow \mathrm{N}(0, \sigma^2(t)\phi_K/f^2(t, Q_\alpha(t))). \tag{21}$$

To construct confidence intervals for $Q_\alpha(t)$ based on (21), one may need to estimate $Q''_\alpha$, which appears highly nontrivial. To circumvent the latter problem, we shall propose a bias-corrected estimate. Let $K^*(u) = 2K(u) - 2^{-1/2}K(u/\sqrt{2})$ and

$$T^*_{\alpha,n}(t) = \sum_{i=1}^{n} \psi_\alpha(X_i - Q_\alpha(t) - Q'_\alpha(t)(i/n - t))K^*_{b_n}(i/n - t).$$

By the argument in (61), we have $\mathbb{E}(T^*_{\alpha,n}(t)) = o(nb_n^3)$ since $\int_{\mathbb{R}} u^2 K^*(u)\,du = 0$. Let

$$\tilde{Q}_{\alpha,b_n}(t) := 2\check{Q}_{\alpha,b_n}(t) - \check{Q}_{\alpha,\sqrt{2}b_n}(t). \tag{22}$$

A similar jackknife-type bias corrected estimate is proposed in Wu and Zhao (2007) for inference of trends in mean nonstationary models. By (20) and following Proposition 6, we have under conditions of Theorem 4 that

$$\sqrt{nb_n}[\tilde{Q}_{\alpha,b_n}(t) - Q_\alpha(t)] = \frac{T^*_{\alpha,n}(t)}{\sqrt{nb_n}} + o_{\mathbb{P}}(1)$$

$$\Rightarrow \mathrm{N}(0, \sigma^2(t)\phi_{K^*}/f^2(t, Q_\alpha(t))). \tag{23}$$

With (22) and (23), to construct confidence intervals for $Q_\alpha(t)$, we can just use $\tilde{Q}_{\alpha,b_n}(t)$ and do not need to estimate the bias term $b_n^2\mu_2 Q''_\alpha(t)/2$ in (21).

3.4. *Estimation of the density and variance functions.* To construct confidence intervals based on Theorems 1 and 4, we should deal with the nontrivial problem of estimating the long-run variance function $\sigma^2(t)$ and the density $f(t, Q_\alpha(t))$. By the local stationarity property, we can estimate them by using observations $X_i$ for which $i/n$ is close to $t$.

For $t \in (0, 1)$, let $s_n(t) = \max(\lfloor nt - nb_n \rfloor, 1)$, $l_n(t) = \min(\lfloor nt + nb_n \rfloor, n)$ and

$$\mathcal{N}_n(t) = \{i \in \mathbb{N} : s_n(t) \leq i \leq l_n(t)\}. \tag{24}$$

Let $Z_{i,\alpha} = \psi_\alpha(X_i - \hat{Q}_\alpha(i/n))$. For a sequence $m_n$ with $m_n \to \infty$ and $nb_n/m_n \to \infty$, let

$$\hat{\sigma}^2(t) = \frac{m_n}{|\mathcal{N}_n(t)| - m_n + 1} \sum_{j=s_n(t)}^{l_n(t)-m_n+1} \left(\frac{1}{m_n}\sum_{i=j}^{j+m_n-1} Z_{i,\alpha} - \bar{Z}_n(t)\right)^2, \tag{25}$$

where $\bar{Z}_n(t) = \sum_{i \in \mathcal{N}_n(t)} Z_{i,\alpha}/|\mathcal{N}_n(t)|$ and $|\mathcal{N}_n(t)| = l_n(t) - s_n(t) + 1$ is the cardinality of $\mathcal{N}_n(t)$. Estimator (25) is a localized version of the popular

block estimate of the long-run variance. For stationary processes, properties of long-run variance estimates have been extensively studied; see, for example, Politis, Romano and Wolf (1999).

For $f(t, Q_\alpha(t))$, we propose the following version of the kernel density estimator

$$\hat{f}(t, Q_\alpha(t)) = \frac{1}{|\mathcal{N}_n(t)| h_n} \sum_{i \in \mathcal{N}_n(t)} K_{h_n}^{\#}(\hat{Q}_\alpha(t) - X_i), \tag{26}$$

where $K^{\#} \in \mathcal{K}$ is a kernel and $h_n$ is the bandwidth satisfying $h_n \to 0$ and $nb_n \times h_n \to \infty$. Theorem 5 asserts consistency of both estimators. Its proof is sketched in Section 6.

THEOREM 5. *Under conditions of Theorem 2, for any $t \in (0, 1)$, $\hat{\sigma}^2(t)$ and $\hat{f}(t, Q_\alpha(t))$ are weakly consistent estimators for $\sigma^2(t)$ and $f(t, Q_\alpha(t))$, respectively.*

Choosing optimal smoothing parameters $m_n$ and $h_n$ is highly nontrivial as it involves the issue of nonstationarity as well as dependence. As a rule of thumb, we recommend using $m_n = \lambda_* |\mathcal{N}_n(t)|^{1/3}$, where $\lambda_*$ can be obtained based on Song (1996), and $h_n = c_* |\mathcal{N}_n(t)|^{-1/5}$, where $c_*$ is chosen by using the bandwidth selector of Sheather and Jones (1991). Those smoothing parameter selectors perform reasonably well for stationary and short-range dependent processes. Since locally $(X_i)$ can be well approximated by a stationary process as the window size $b_n \to 0$, it is expected that those selectors are also applicable under our setting.

**4. Examples.** In this section we present examples of locally stationary linear and nonlinear time series. For such processes, the stability of $F_k(\cdot, \cdot, \Upsilon_i)$ as well as the local stationarity conditions can be verified. Hence, results in Section 3 are applicable.

4.1. *Nonstationary linear processes.* Let $\varepsilon_i$ be i.i.d. random variables; let $a_j(\cdot), j = 0, 1, \ldots,$ be $\mathcal{C}^1[0, 1]$ functions such that

$$G(t, \Upsilon_i) = \sum_{j=0}^{\infty} a_j(t) \varepsilon_{i-j} \tag{27}$$

is well defined for all $t \in [0, 1]$. Model (27) was considered in Dahlhaus (1997), where his primary interest is to estimate time-varying spectra.

Propositions 2 and 3 are for checking the stability condition and the stochastic Lipschitz continuity condition (A3), respectively. They are closely related the example in Section 2.1 in Draghicescu, Guillas and Wu (2008).

Here for the sake of completeness we provide rigorous and explicit statements. Let $F_\varepsilon$ be the distribution function of $\varepsilon_i$ and $f_\varepsilon$ be its density. For $k = 0, 1, \ldots$, denote by $g^{(k)}(\cdot)$ the $k$th derivative of $g(\cdot) : \mathbb{R} \to \mathbb{R}$. Here $g^{(0)} = g$.

PROPOSITION 2. *Assume that $\varepsilon_0 \in \mathcal{L}^q$, $q > 0$, and its density $f_\varepsilon(\cdot)$ satisfies*

$$\sup_x |f_\varepsilon^{(k)}(x)| < C_0, \qquad k = 0, 1, 2, 3, \tag{28}$$

*for some $C_0 < \infty$. Further assume $\min_{t \in [0,1]} |a_0(t)| > 0$. Let $q' = \min(2, q)$. Then $\delta_{F_k}(i) = O(\sup_{t \in [0,1]} |a_i(t)|^{q'/2})$. Hence $\{F_k(\cdot, \cdot, \Upsilon_i)\}_{i \in \mathbb{Z}}$ is stable if $\sum_{j=0}^\infty \sup_{t \in [0,1]} |a_j(t)|^{q'/2} < \infty$.*

PROOF. We first assume $a_0(t) \equiv 1$. Let $\underline{G}(t, \Upsilon_{i-1}) = G(t, \Upsilon_i) - \varepsilon_i$ and $\underline{G}(t, \Upsilon_{i-1}^*) = \underline{G}(t, \Upsilon_{i-1}) - a_i(t)\varepsilon_0 + a_i(t)\varepsilon_0'$. Using $\min(1, |x|) \leq |x|^{q'/2}$, we have

$$\begin{aligned}
\delta_{F_k}(i) &= \sup_{t \in [0,1]} \sup_{u \in \mathbb{R}} \|F_\varepsilon^{(k)}(u - \underline{G}(t, \Upsilon_{i-1})) - F_\varepsilon^{(k)}(u - \underline{G}(t, \Upsilon_{i-1}^*))\| \\
&\leq \sup_{t \in [0,1]} \min\{2C_0, C_0 \|a_i(t)\varepsilon_0 - a_i(t)\varepsilon_0'\|\} \\
&\leq \sup_{t \in [0,1]} 2C_0 \||a_i(t)\varepsilon_0 - a_i(t)\varepsilon_0'|^{q'/2}\| = O\left(\sup_{t \in [0,1]} |a_i(t)|^{q'/2}\right).
\end{aligned}$$

The case in which $a_0(t) \not\equiv 1$ can be similarly dealt with by using the argument in Section 2.1 in Draghicescu, Guillas and Wu (2008). □

PROPOSITION 3. *Assume that $\varepsilon_0 \in \mathcal{L}^q$, $q > 0$. Then condition* (A3) *holds with this $q$ if*

$$\sum_{j=0}^\infty \sup_{t \in [0,1]} |a_j'(t)|^{\min(2,q)} < \infty. \tag{29}$$

Note that $\zeta_i(t) - \zeta_i(s) = \sum_{j=0}^\infty [a_j(t) - a_j(s)]\varepsilon_{i-j}$. Using the Lebesgue Dominant Convergence Theorem and similar arguments as those in the proof of Proposition 2, Proposition 3 easily follows. Details are omitted. Note that both propositions allow heavy-tailed processes $\{X_i\}$.

4.2. *Nonstationary nonlinear time series.* Let $\varepsilon_i$ be i.i.d. Many stationary nonlinear time series models are of the form

$$Z_i = R(Z_{i-1}, \varepsilon_i), \tag{30}$$

where $R$ is a measurable function [Wu and Shao (2004)]. A natural extension of (30) to the locally stationary setting is to incorporate the time index $t$ via

$$\zeta_i(t) = R(t, \zeta_{i-1}(t), \varepsilon_i), \qquad 0 \le t \le 1. \tag{31}$$

Theorem 6 asserts that, under suitable conditions, (31) has a unique stationary solution of the form $\zeta_i(t) = G(t, \Upsilon_i)$. Then one can have a nonstationary process $X_{i,n} = G(i/n, \Upsilon_i)$.

THEOREM 6. *Let $q > 0$. Assume that, for some $x_0$, $\sup_{t \in [0,1]} \|R(t, x_0, \varepsilon_i)\|_q < \infty$, and*

$$\chi := \sup_{t \in [0,1]} L(t) < 1, \qquad \text{where } L(t) = \sup_{x \ne y} \frac{\|R(t, x, \varepsilon_0) - R(t, y, \varepsilon_0)\|_q}{|x - y|}. \tag{32}$$

*Then for any $t \in [0,1]$, (31) admits a unique stationary solution, and iterations of (31) lead to $\zeta_i(t) = G(t, \Upsilon_i)$. Furthermore, we have*

$$\sup_{t \in [0,1]} \|G(t, \Upsilon_i) - G(t, \Upsilon_i^*)\|_q = O(\chi^i). \tag{33}$$

Theorem 6 follows from the argument of Theorem 2 in Wu and Shao (2004). We omit the proof since there is no essential extra difficulties involved. By the Markovian structure of $(\zeta_i(t))$, we can write $F(t, x, \Upsilon_i) = F(t, x, \zeta_i(t))$. Let $\dot{F}_k(t, x, u) = \partial F_k(t, x, u)/\partial u$ if the latter exists. The following corollary is immediate [cf. Proposition 5 in Wu (2007a)].

COROLLARY 2. *Assume that conditions of Theorem 6 hold. Further assume*

$$\sup_{t \in [0,1]} \sup_{x,u} |\dot{F}_k(t, x, u)| + \sup_{t \in [0,1]} \sup_{x,u} |F_k(t, x, u)| < \infty, \qquad k = 0, 1, 2, 3. \tag{34}$$

*Then we have $\delta_{F_k}(i) = O(\chi_1^i)$ for some $\chi_1 \in (0,1)$. Hence $\{F_k(\cdot, \cdot, \Upsilon_i)\}_{i \in \mathbb{Z}}$ is stable.*

PROPOSITION 4. *Under conditions of Theorem 6, we have* (A3) *provided that*

$$\sup_{t \in [0,1]} \|M(G(t, \Upsilon_0))\|_q < \infty,$$

$$\text{where } M(x) = \sup_{0 \le t < s \le 1} \frac{\|R(t, x, \varepsilon_0) - R(s, x, \varepsilon_0)\|_q}{|t - s|}. \tag{35}$$

PROOF. We only consider the case $q \geq 1$. The case $q < 1$ similarly follows. From the proof of Theorem 2 in Wu and Shao (2004), we have $\sup_{t \in [0,1]} \|G(t, \Upsilon_0)\|_q < \infty$. Let $C_1 = \sup_{t \in [0,1]} \|M(G(t, \Upsilon_0))\|_q$ and $\chi = \sup_{t \in [0,1]} L(t) < 1$. For $t, s \in [0,1]$, by (32),

$$\|R(s, \zeta_{i-1}(t), \varepsilon_i) - R(s, \zeta_{i-1}(s), \varepsilon_i)\|_q \leq \chi \|\zeta_{i-1}(t) - \zeta_{i-1}(s)\|_q.$$

Hence, by (35), $\|\zeta_i(t) - \zeta_i(s)\|_q \leq C_1|t-s| + \chi \|\zeta_i(t) - \zeta_i(s)\|_q$. So (A3) follows. □

EXAMPLE 1 (Time varying threshold autoregressive (TVTAR) models). Let $\varepsilon_i \in \mathcal{L}^q$, $q > 0$, be i.i.d. with distribution function $F_\varepsilon$ and density $f_\varepsilon = F'_\varepsilon$. Consider the model

$$(36) \quad \zeta_i(t) = a(t)[\zeta_{i-1}(t)]^+ + b(t)[-\zeta_{i-1}(t)]^+ + \varepsilon_i, \qquad 0 \leq t \leq 1,$$

where $a(t), b(t) \in \mathcal{C}^1[0,1]$. Then Theorem 6 is applicable if $\sup_{t \in [0,1]}[|a(t)| + |b(t)|] < 1$. Since $F(t, x, u) = F_\varepsilon(x - a(t)u^+ - b(t)(-u)^+)$, we have (34) if there exits a constant $C < \infty$, such that $\sup_u |f_\varepsilon^{(k)}(u)| < C$, $k = 0, \ldots, 3$. Since $a(t), b(t) \in \mathcal{C}^1[0,1]$, (35) holds. Note that, similarly as the linear process case, $(X_i)$ can be heavy-tailed.

**5. Data analysis.** Climate change has received enormous attention. The understanding of changes in trends, variability and extremes of climate is of great importance. See Contribution of Working Group I to the fourth Assessment Report of the Intergovernmental Panel on Climate Change Solomon et al. (2007) (WGI07 hereafter) for a comprehensive discussion.

The local linear quantile smoothing technique is a useful approach to studying changes in trends and variability of nonstationary time series by examining quantile curves, for example, the median and IQR curves. As pointed out in WGI07, extremes are infrequent events at the high and low ends of the range of values of a particular variable. Hence, upper and lower quantile curves such as $Q_{0.95}(t)$ and $Q_{0.05}(t)$ are natural tools for illustrating changes in climate extremes.

Robust tools such as quantile methods have already been applied in the climatology community. See, for instance, Hannachi (2006). However, previous literature is mostly for stationary or integrated linear time series models. Additionally, Gaussian or other specific distributional assumptions were typically made. These assumptions are vulnerable to model misspecifications, which may lead to erroneous conclusions. Here we shall apply our method to climate data sets and compare our findings with previous ones.

5.1. *Global temperature data.* Global warming is one of the most important issues in climate change and it is closely related to the changes of other climate variables (sea level, for instance). Global temperature series have been extensively studied in the statistics community; see, for example, Bloomfield and Nychka (1992), Vogelsang (1998), Wu, Woodroofe and Mentz (2001) and Wu and Zhao (2007), among others. Most of the above studies focused on estimation and inference of trends (mean functions). Less statistical research has been done with respect to changes in variability and extremes of global temperatures.

Here we consider the series complied by Jones et al.; see http://cdiac.esd.ornl.gov/ftp/trends/temp/jonescru/. It contains global monthly temperature anomalies from 1856 to 2005, relative to the 1961–1990 mean; see left panel of Figure 2 for the data. We study the median, 95% and 5% extreme quantile and IQR curves. Jackknife bias correction and two-stage smoothing are applied. We tried both the static bandwidth (15) and the time-varying $b_n(\alpha, t)$. Both choices of bandwidths yield quite similar results. Hence here we present our results by using the static bandwidth. For the 0.05th, 0.25th, 0.5th, 0.75th and 0.95th quantile curves, the bandwidths are chosen as 0.083, 0.077, 0.075, 0.077 and 0.089, respectively. The bandwidth for the second-stage smoothing step is 0.04. The estimated curves and corresponding asymptotic 95% point-wise confidence bands are shown in Figure 1.

From the estimated median curve, we see a clear warming trend. Interestingly, the increasing of temperature is not homogenous, with a strong faster warming trend after 1975. Wu and Zhao (2007) rejected the hypothesis of linear increasing of trend. Our estimated median curve is consistent with the mean temperature curve in Figure TS.6 in WGI07. Similar trends of significant nonhomogenous increase are found in both 5% and 95% quantile curves of monthly temperature (although 95% quantile curve seems to slightly decrease after 1995). As summarized in WGI07,

> "Changes in extremes of temperature are consistent with warming. Observations show widespread reductions in the number of frost days in mid-latitude regions, increases in the number of warm extremes (warmest 10% of days or nights) and a reduction in the number of daily cold extremes (coldest 10% of days or nights)."

Note that the magnitude of increase of the 5% quantile is approximately twice that of the 95% quantile, indicating reduction in the temperature variability. This result is consistent with the statement in *Climate Change 2001: The Scientific Basis* [IPCC (2001), WGI01 hereafter], that while increases in the frequency of warm days have been observed, decreases in the number of cool nights have been stronger. Reduction in the variability can be verified by the IQR curve. Compared with the width of the confidence band, the reduction of variability is not as significant as the increasing of median and extreme quantiles. In Chapter 2 of WGI01 several climatological papers were cited to confirm that annual and monthly variation of global temperature had decreased. To summarize, although based on different statistical models and assumptions, results of our method are consistent with the previous empirical findings. Our results suggest that global monthly temperature distribution in the last 150 years has the tendency of shifting to the right with shrinking variability.

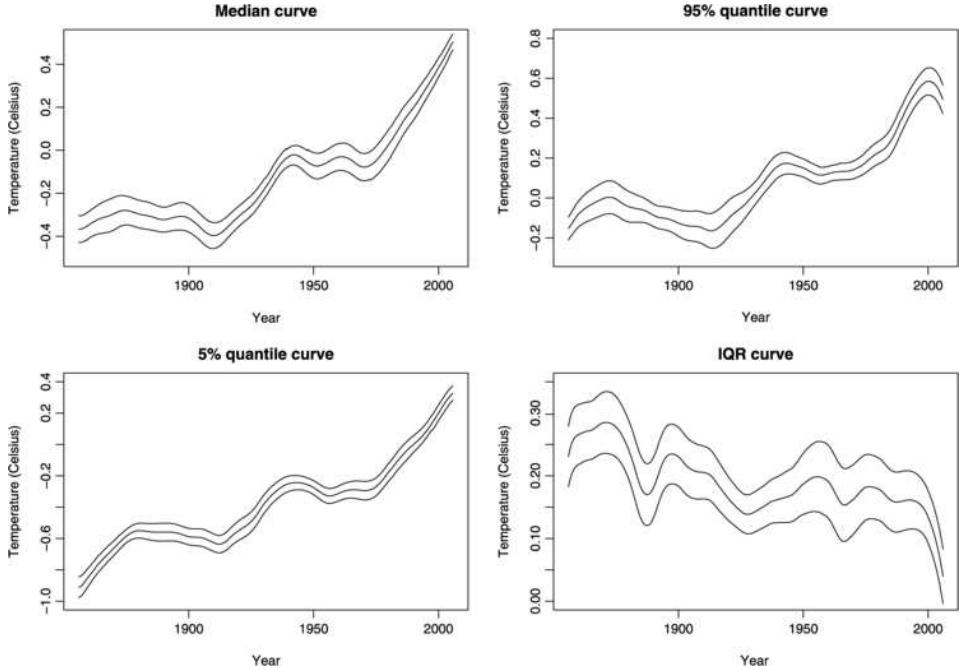

FIG. 1. *Trends in median, extreme quantiles and IQR of the global monthly temperature series, with* 95% *point-wise confidence bands.*

5.2. *U.S. precipitation data.* This data set contains spatially averaged monthly total precipitation in the land surface of United States from January 1895 to April 2007. It is available at National Climate Data Center's website at http://www1. ncdc.noaa.gov/pub/data/cirs/drd964x.pcpst.txt. One can refer to http://www1.ncdc.noaa.gov/pub/data/cirs/state.README for a detailed description.

Average seasonal trend was first removed from the data and we focus on inference of the anomalies; the right panel of Figure 2 summarizes the anomalies. Jackknife bias correction and two-stage smoothing are applied. We also choose the static bandwidth (15) since the time-varying bandwidths provide similar estimated curves. For the 0.05th, 0.25th, 0.5th, 0.75th and 0.95th quantile curves, the bandwidths are chosen as 0.114, 0.096, 0.096, 0.096 and 0.115, respectively. The bandwidth for the second-stage smoothing step is chosen as 0.05. The estimated median, extreme quantile and IQR curves and corresponding asymptotic 95% point-wise confidence bands are shown in Figure 3.

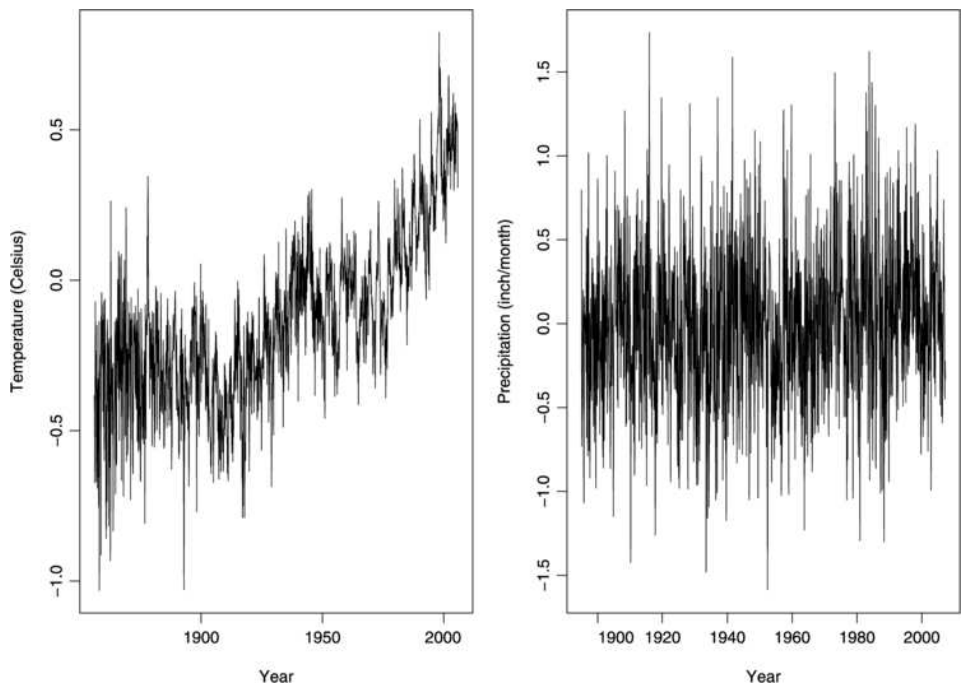

FIG. 2. *Time series plots of the data examples. Left panel: global warming data. Right panel: U.S. precipitation data.*

In theory, increased temperature is likely to produce heavier precipitation. As stated in WGI01,

> "Increasing global surface temperatures are very likely to lead to changes in precipitation and atmospheric moisture, because of changes in atmospheric circulation, a more active hydrological cycle, and increases in the water holding capacity throughout the atmosphere."

However, the estimated median curve shows that even though there is a slight increase of the monthly precipitation, noticeably after 1960, the trend is not significant with respect to its stochastic uncertainty. WGI01 suggested to compare changes in many of the moisture-related variables, such as streamflow and soil moisture, with precipitation to help validate long-term precipitation trends and claimed that the precipitation change was not spatially uniform. On the other hand, one factor which could be considered is the so-called "global dimming" effect [see, e.g., Stanhill and Cohen (2001), Roderick and Farquhar (2002), Liepert et al. (2004), among others], which reduces the sun radiation and water evaporation on the surface and therefore reduces precipitation. Hence future research taking into account of those factors may help confirm (or reject) the heavier precipitation hypothesis. Both WGI07 and WGI01 claimed that there was an amplified increase of heavy precipitation, even in regions where the mean precipitation decreased. Our estimated 95% quantile curve does show a stronger increase than the mean, but still not very impressive considering the width of the confidence band. The 5% quantile and IQR curves also show no significant increase or decrease pattern. To summarize, the U.S. monthly precipitation data does not show significant increasing trend of precipitation median and heavy precipitation, as reported in the literature and a more comprehensive study including more factors may be helpful.

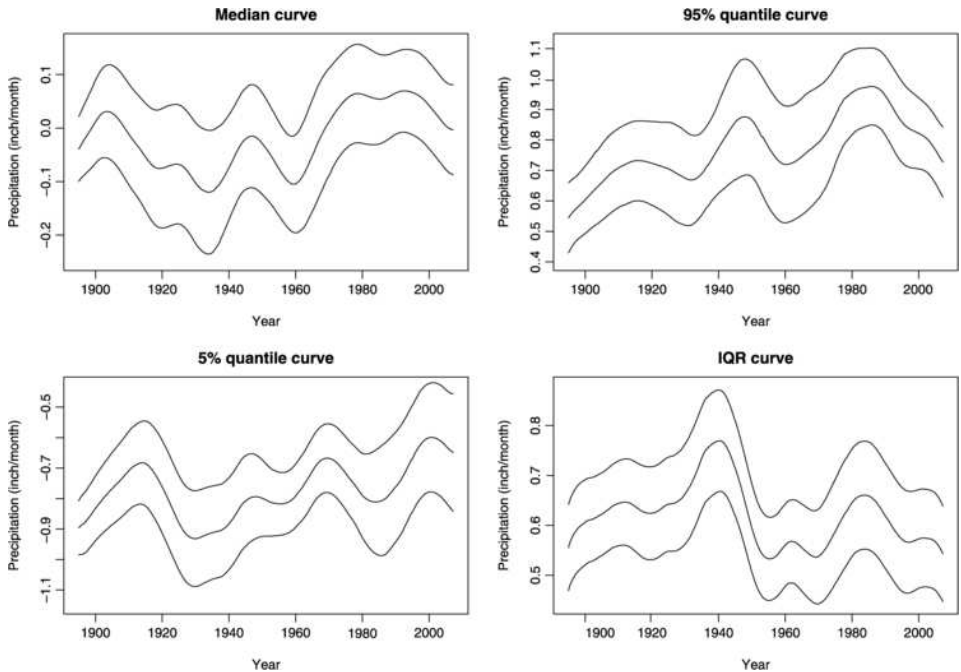

Fig. 3. *Trends in median, extreme quantiles and IQR of monthly U.S. precipitation, with 95% point-wise confidence bands.*

**6. Proof of results in Section 3.** Unless otherwise specified, we will only prove results in Section 3 for $\alpha = 1/2$, since results for other quantiles similarly follow. We shall also omit the subscript $\alpha$ in the notation if no confusion will be caused. Here we list some notation that we will use in this section. Let $Y_i(t) = X_i - Q(t) - (i/n - t)Q'(t)$, $\check{Y}_i(t) = G(t, \Upsilon_i) - Q(t)$, $\check{e}_i(t) = \psi(\check{Y}_i(t))$ and $\check{S}_n(t) = \sum_{i=1}^n \check{e}_i(t)\mathbf{v}_{i,n}(t)$, where $\mathbf{v}_{i,n}(t) = K_{b_n}(i/n - t)\mathbf{z}_{i,n}(t)$. Recall $\mathbf{z}_{i,n}(t) = (1, (i/n - t)/b_n)^\top$.

For $k \in \mathbb{Z}$, define the projection operator $\mathcal{P}_k$ on $\mathcal{L}^1$, by $\mathcal{P}_k \cdot = \mathbb{E}(\cdot|\Upsilon_k) - \mathbb{E}(\cdot|\Upsilon_{k-1})$. Recall $J(t, x, \Upsilon_i) = I\{G(t, \Upsilon_i) \leq x\}$. Let $i \geq 1$. Then $\mathbb{E}(J(t, x, \Upsilon_i)|\Upsilon_0) = \mathbb{E}(F(t, x, \Upsilon_{i-1})|\Upsilon_0)$ and $\mathbb{E}(J(t, x, \Upsilon_i)|\Upsilon_{-1}) = \mathbb{E}(F(t, x, \Upsilon_{i-1})|\Upsilon_{-1})$. Since $\varepsilon_0', \varepsilon_j, j \in \mathbb{Z}$, are independent, we have $\mathbb{E}(F(t, x, \Upsilon_{i-1})|\Upsilon_{-1}) = \mathbb{E}(F(t, x, \Upsilon_{i-1}^*)|\Upsilon_{-1})$ and $\mathbb{E}(F(t, x, \Upsilon_{i-1}^*)|\Upsilon_{-1}) = F(t, x, \Upsilon_{i-1}^*)|\Upsilon_0)$. Hence

$$
(37) \quad \|\mathcal{P}_0 J(t, x, \Upsilon_i)\| = \|\mathbb{E}(F(t, x, \Upsilon_{i-1}) - F(t, x, \Upsilon_{i-1}^*)|\Upsilon_0)\|
$$
$$
\leq \|F(t, x, \Upsilon_{i-1}) - F(t, x, \Upsilon_{i-1}^*)\|.
$$

Since $\{F(\cdot, \cdot, \Upsilon_i)\}_{i \in \mathbb{Z}}$ is stable, by (5), $\sigma(t)$ in (7) is finite.

LEMMA 1. *Let $b_n \to 0$; let $(\alpha_{ni})_{i=1}^n$ be a triangular array of real numbers such that*

$$
(38) \quad A_n(t) := \sum_{i=1}^n \alpha_{ni}^2 K_{b_n}^2(i/n - t) \leq C,
$$

*and $\varpi_n := \max_{i \in \mathcal{N}_n(t)} |\alpha_{ni}| \to 0$. Assume (A1), (A2) and that $\{F(\cdot, \cdot, \Upsilon_i)\}_{i \in \mathbb{Z}}$ is stable. Let $\eta_{ni}(t) = [\rho(Y_i(t) - \alpha_{ni}) - \rho(Y_i(t)) + \alpha_{ni}\psi(Y_i(t))]K_{b_n}(i/n - t)$. Then $\mathrm{var}[\sum_{i=1}^n \eta_{ni}(t)] \to 0$.*

PROOF. Let $\mathcal{A} = [\inf_{t \in [0,1]} Q(t) - 1, \sup_{t \in [0,1]} Q(t) + 1]$. For every $k \geq 1$, sufficiently small $\epsilon_0$ and sufficiently large $n$, we have by (37) that, for $|\epsilon| \leq \epsilon_0$,

$$
(39) \quad \mathrm{SUP}\|\mathcal{P}_{i-k}[\psi(Y_i(t) - \epsilon) - \psi(Y_i(t))]\| \leq \mathrm{SUP} 2\|\mathcal{P}_{i-k}\psi(Y_i(t) - \epsilon)\|
$$
$$
\leq 2\delta_F(k-1),
$$

where SUP denotes $\sup_{|\epsilon| \leq \epsilon_0} \sup_{i \in \mathcal{N}_n(t)}$. On the other hand, for sufficiently small $\epsilon_0$ and sufficiently large $n$, by condition (A2),

$$
(40) \quad \mathrm{SUP}\|\mathcal{P}_{i-k}[\psi(Y_i(t) - \epsilon) - \psi(Y_i(t))]\|
$$
$$
\leq \mathrm{SUP}\|\psi(Y_i(t) - \epsilon) - \psi(Y_i(t))\|
$$
$$
\leq \left(\epsilon_0 \max_{(t,x) \in [0,1] \times \mathcal{A}} f(t, x)\right)^{1/2} \leq C\epsilon_0^{1/2}.
$$

Let $Z_{k,n}(s) = \sum_{i=1}^n \mathcal{P}_{i-k}[\psi(Y_i(t) - s\alpha_{ni}) - \psi(Y_i(t))]\alpha_{ni}K_{b_n}(i/n - t)$. Since the summands of $Z_{k,n}(s)$ are martingale differences, by (38), (39) and (40), we have, for large $n$,

$$\sup_{0\leq s\leq 1} \|Z_{k,n}(s)\|^2 = \sup_{0\leq s\leq 1} \sum_{i\in\mathcal{N}_n(t)} \|\mathcal{P}_{i-k}[\psi(Y_i(t)-s\alpha_{ni})-\psi(Y_i(t))]\|^2$$

(41)
$$\times \alpha_{ni}^2 K_{b_n}^2(i/n-t)$$

$$\leq C\min\{\delta_F^2(k-1),\varpi_n\}$$

holds for $k\geq 1$. Clearly, $\sup_{0\leq s\leq 1}\|Z_{0,n}(s)\|^2 \leq C\min(1,\varpi_n)$. Since $\{F(\cdot,\cdot,\Upsilon_i)\}$ is stable, $\sum_{k=1}^\infty \delta_F(k-1) < \infty$. By the Lebesgue Dominated Convergence Theorem, we have

$$\left\|\sum_{i=1}^n [\eta_{ni}(t)-\mathbb{E}\eta_{ni}(t)]\right\| = \left\|-\sum_{k=0}^\infty \int_{s=0}^1 Z_{k,n}(s)\,ds\right\| \leq \sum_{k=0}^\infty \int_{s=0}^1 \|Z_{k,n}(s)\|\,ds$$

(42)
$$\leq C\min(1,\varpi_n) + \sum_{k=1}^\infty C\min\{\delta_F(k-1),\varpi_n^{1/2}\}$$

$$\to 0$$

as $n\to\infty$ in view of $\varpi_n \to 0$. □

LEMMA 2.  *Assume* (A1), (A2) *and* $b_n\to 0$. *Then for all large* $n$,

(43)
$$\sup_{i\in\mathcal{N}_n(t)} |\mathbb{E}[\rho(Y_i(t)+a)-\rho(Y_i(t))-a\psi(Y_i(t))] - f(t,Q(t))a^2/2|$$
$$= o(a^2) \quad \text{as } a\to 0.$$

PROOF. Let $T\in\mathbb{N}$ and $\tau = a/T$. Since $\rho(x)$ is convex, for all $j=1,\ldots,T$, we have

(44)
$$\tau[\psi(Y_i(t)+j\tau-\tau)-\psi(Y_i(t))]$$
$$\leq \rho(Y_i(t)+j\tau) - \rho(Y_i(t)+j\tau-\tau) - \psi(Y_i(t))\tau$$
$$\leq \tau[\psi(Y_i(t)+j\tau)-\psi(Y_i(t))].$$

Taking expectation, using (A2), summing over $j=1,\ldots,T$, and then letting $T\to\infty$, we have (43). □

Before we proceed to the next proposition, let us first observe that the original estimation problem (3) is equivalent to the following transformed one:

(45)
$$\hat{\vartheta}(t) := (\hat{\vartheta}_1(t),\hat{\vartheta}_2(t))^\top$$
$$= \arg\min_\vartheta \sum_{i=1}^n \rho(Y_i(t)-\vartheta^\top \mathbf{u}_{i,n}(t))K_{b_n}(i/n-t),$$

where

(46) $\quad \mathbf{u}_{i,n}(t) := (u_{i,n,1}(t), u_{i,n,2}(t))^\top = \Sigma_n(t)^{-1/2}(1, (i/n - t)/b_n)^\top,$

and $\Sigma_n(t) = \sum_{i=1}^n (1, (i/n - t)/b_n)^\top (1, (i/n - t)/b_n) K_{b_n}(i/n - t)$. Since $K \in \mathcal{K}$, we have $\Sigma_n(t) = nb_n \mu_{\mathbf{K}} + O(1)$ uniformly on any closed interval of $(0, 1)$.

The estimators $(\hat{Q}(t), \hat{Q}'(t))^\top$ of (3) and $\hat{\vartheta}(t)$ are related via

(47)
$$\hat{\vartheta}_1(t) = \Sigma_n(t)^{1/2}(\hat{Q}(t) - Q(t)) \quad \text{and}$$
$$\hat{\vartheta}_2(t) = \Sigma_n(t)^{1/2} b_n (\hat{Q}'(t) - Q'(t)).$$

Recall $\hat{\theta}(t) = (\hat{Q}(t) - Q(t), b_n(\hat{Q}'(t) - Q'(t)))^\top$ and (16) for $S_n(t)$.

PROPOSITION 5. *Assume* (A1), (A2), $b_n \to 0$, $nb_n \to \infty$, $f(t, Q(t)) > 0$ *and that* $\{F(\cdot, \cdot, \Upsilon_i)\}_{i \in \mathbb{Z}}$ *is stable. Then*

(48) $\quad f(t, Q(t))\Sigma_n(t)\hat{\theta}(t) - S_n(t) = o_\mathbb{P}((nb_n)^{1/2}).$

PROOF. For a fixed vector $\vartheta \in \mathbb{R}^2$, let $\alpha_{ni} = \vartheta^\top \mathbf{u}_{i,n}(t)$. Then $\max_{i \in \mathcal{N}_n(t)} |\alpha_{ni}| \to 0$. Note that $\sum_{i=1}^n \mathbf{u}_{i,n}(t)\mathbf{u}_{i,n}^\top(t) K_{b_n}(i/n - t) = Id_2$, where $Id_2$ is the $2 \times 2$ identity matrix. Let $\varpi_i(t) = \rho(Y_i(t) - \vartheta^\top \mathbf{u}_{i,n}(t)) - \rho(Y_i(t)) + \vartheta^\top \mathbf{u}_{i,n}(t)\psi(Y_i(t))$. By Lemma 2, we have

(49)
$$\sum_{i=1}^n \mathbb{E}[\varpi_i(t)] K_{b_n}(i/n - t) = \sum_{i \in \mathcal{N}_n(t)} f(t, Q(t)) \frac{(\vartheta^\top \mathbf{u}_{i,n}(t))^2}{2}$$
$$\times (1 + o(1)) K_{b_n}(i/n - t)$$
$$= \frac{f(t, Q(t))|\vartheta|^2}{2} + o(|\vartheta|^2).$$

Since $K \in \mathcal{K}$, by simple calculations, (38) holds. Thus by Lemma 1 and (49), we have

(50) $\quad \sum_{i=1}^n \varpi_i(t) K_{b_n}(i/n - t) \to \frac{f(t, Q(t))|\vartheta|^2}{2} \quad \text{in probability}.$

Note that both sides of (50) are convex functions of $\vartheta$. By the Convexity Lemma in Pollard (1991), page 187, we have that the above convergence holds uniformly for $|\vartheta| \le c$, where $c > 0$ is a finite constant. Now the relation (48) easily follows from properties of convex functions; see, for example, the proofs of Theorem 1 in Wu (2007a), Theorems 2.2 and 2.4 in Bai, Rao and Wu (1992) or Theorem 1 in Pollard (1991). Details are omitted. $\square$

LEMMA 3. *Let* $\check{V}_n(t) = \operatorname{cov}(\check{S}_n(t))$ *and* $\nu_{\mathbf{K}} = \operatorname{diag}(\int_{-1}^1 K^2(u)\,du, \int_{-1}^1 u^2 \times K^2(u)\,du)$. *Suppose that* $\{F(\cdot, \cdot, \Upsilon_i)\}_{i \in \mathbb{Z}}$ *is stable*, $b_n \to 0$ *and* $nb_n \to \infty$. *Then*

(51) $\quad (nb_n)^{-1}\check{V}_n(t) \to \sigma^2(t)\nu_{\mathbf{K}}.$

PROOF. For simplicity we write $\check{e}_i$ etc for $\check{e}_i(t)$. Clearly $(\check{e}_i)_{i\in\mathbb{Z}}$ is stationary and $\mathbb{E}(\check{e}_i) = 0$. Let $D_i = \sum_{j=i}^{\infty} \mathcal{P}_i \check{e}_j$ and $R_i = \check{e}_i - D_i$. Recall (24) for $\mathcal{N}_n(t)$. Let $s_n = s_n(t)$, $l_n = l_n(t)$, $A_i = \sum_{j=s_n}^{s_n+i}(\check{e}_j - D_j)$ and $B_n = \sum_{i\in\mathcal{N}_n(t)} D_i \mathbf{v}_{i,n}$. By the summation by parts formula,

$$\|\check{S}_n(t) - B_n\| = \left\|\sum_{i\in\mathcal{N}_n(t)} R_i \mathbf{v}_{i,n}\right\|$$

$$\leq \sum_{i=0}^{l_n-s_n-1} \|A_i\| |\mathbf{v}_{i,n} - \mathbf{v}_{i+1,n}| + \|A_{l_n-s_n}\| |\mathbf{v}_{l_n-s_n,n}|.$$

Let $\Omega_m = \sum_{k=m}^{\infty} \delta_F(k-1)$. By (3) of Theorem 1 in Wu (2007b), we have

$$\omega_n := \max_{0\leq i\leq l_n-s_n} \|A_i\|^2 \leq C \sum_{j=1}^{l_n-s_n} \Omega_j^2 = o(l_n - s_n) = o(nb_n)$$

since $\Omega_m \to 0$ as $m \to \infty$. Hence

$$(52)\quad \|\check{S}_n(t) - B_n\| \leq \omega_n^{1/2} \left(\sum_{i=0}^{l_n-s_n-1} |\mathbf{v}_{i,n} - \mathbf{v}_{i+1,n}| + |\mathbf{v}_{l_n-s_n,n}|\right) = o((nb_n)^{1/2}).$$

Note that $D_i$ are martingale differences. By orthogonality, since $K \in \mathcal{K}$, it is easily seen that $(nb_n)^{-1}\mathbb{E}(B_n B_n^\top) \to \sigma^2(t)\nu_\mathbf{K}$. Hence, by (52), we have (51). □

LEMMA 4. *Under conditions of Lemma 3 and assuming $\sigma(t) > 0$, we have*

$$(53)\quad (nb_n)^{-1/2}\check{S}_n(t) \Rightarrow \mathrm{N}(0, \sigma^2(t)\nu_\mathbf{K}).$$

PROOF. For $m \in \mathbb{N}$, let $\tilde{e}_{m,i} = \mathbb{E}(\check{e}_i | \Upsilon_{i-m,i})$, where, for $j \leq k$, $\Upsilon_{j,k} = (\varepsilon_j, \varepsilon_{j+1}, \ldots, \varepsilon_k)$. Here we also write $\check{e}_i$ etc for $\check{e}_i(t)$. Then $(\tilde{e}_{m,i})_{i=1}^n$ is $m$-dependent. Elementary calculations show that $\mathcal{P}_{i-k}\tilde{e}_{m,i} = \mathbb{E}(\mathcal{P}_{i-k}\check{e}_i | \Upsilon_{i-m,i})$. Thus, we have for sufficiently large $n$,

$$\sup_{i\in\mathcal{N}_n(t)} \|\mathcal{P}_{i-k}(\check{e}_i - \tilde{e}_{m,i})\| \leq 2 \sup_{i\in\mathcal{N}_n(t)} \|\mathcal{P}_{i-k}\check{e}_i\| \leq 2\delta_F(k-1).$$

On the other hand, by the stationarity of $(\check{e}_i)$,

$$\sup_{i\in\mathcal{N}_n(t)} \|\mathcal{P}_{i-k}(\check{e}_i - \tilde{e}_{m,i})\| \leq \sup_{i\in\mathcal{N}_n(t)} \|\check{e}_i - \tilde{e}_{m,i}\| = \|\check{e}_0 - \tilde{e}_{m,0}\| := \tau_J(m).$$

Clearly $\lim_{m\to\infty} \tau_J(m) = 0$. By the orthogonality,

$$\left\|\sum_{i=1}^n \mathcal{P}_{i-k}(\check{e}_i - \tilde{e}_{m,i})\mathbf{v}_{i,n}\right\|^2 = \sum_{i=1}^n \left\|\mathcal{P}_{i-k}(\check{e}_i - \tilde{e}_{m,i})\right\|^2 |\mathbf{v}_{i,n}|^2$$

$$\text{(54)} \qquad \leq \min\{\tau_J^2(m), \delta_F^2(k-1)\} \sum_{i=1}^n |\mathbf{v}_{i,n}|^2$$

$$\leq Cnb_n \min\{\tau_J^2(m), \delta_F^2(k-1)\}.$$

Write $\tilde{S}_{m,n}(t) = \sum_{i=1}^n \tilde{e}_{m,i} \mathbf{v}_{i,n}$. Then $\check{S}_n(t) - \tilde{S}_{m,n}(t) = \sum_{k=0}^\infty \sum_{i=1}^n \mathcal{P}_{i-k}(\check{e}_i - \tilde{e}_{m,i}) \mathbf{v}_{i,n}$. Let LIM denote $\lim_{m\to\infty} \limsup_{n\to\infty}$. Since $\{F(\cdot,\cdot,\Upsilon_i)\}_{i\in\mathbb{Z}}$ is stable, as (42), by the Lebesgue Dominated Convergence Theorem,

$$\text{LIM} \frac{\|\check{S}_n(t) - \tilde{S}_{m,n}(t)\|}{(nb_n)^{1/2}}$$

$$= \text{LIM} \frac{\|\sum_{k=0}^\infty \sum_{i=1}^n \mathcal{P}_{i-k}(\check{e}_i - \tilde{e}_{m,i}) \mathbf{v}_{i,n}\|}{(nb_n)^{1/2}}$$

$$\text{(55)} \qquad \leq \text{LIM}(nb_n)^{-1/2} \sum_{k=0}^\infty \left\|\sum_{i=1}^n \mathcal{P}_{i-k}(\check{e}_i - \tilde{e}_{m,i}) \mathbf{v}_{i,n}\right\|$$

$$\leq \text{LIM}\left[C \min(1, \tau_J(m)) + \sum_{k=1}^\infty C \min\{\tau_J(m), 2\delta_F(k-1)\}\right]$$

$$\to 0.$$

For a symmetric matrix $A$ let $\lambda(A)$ be its smallest eigenvalue. Let $\tilde{V}_{m,n}(t) = \text{cov}(\tilde{S}_{m,n}(t))$. By (55) and Lemma 3, $\text{LIM}|\sigma^2(t)\nu_\mathbf{K} - (nb_n)^{-1}\tilde{V}_{m,n}(t)| = 0$. Note that $\sigma^2(t)\nu_\mathbf{K}$ is nonsingular. Thus there exists an $m_0 \in \mathbb{N}$ and $c_0 > 0$ such that for $m \geq m_0$ and for all large $n$, $(nb_n)^{-1}\lambda(\tilde{V}_{m,n}(t)) > c_0$. Note that $(\tilde{e}_{m,i})$ are $m$-dependent. Applying the Cramer–Wold device and the central limit theorem for $m$-dependent random variables [see Hoeffding and Robbins (1948)], we have $\tilde{V}_{m,n}^{-1/2}(t)[\tilde{S}_{m,n}(t) - \mathbb{E}\tilde{S}_{m,n}(t)] \Rightarrow \text{N}(0, Id_2)$ for fixed $m \geq m_0$. Thus by (55), we obtain (53). □

PROPOSITION 6. *Under conditions of Theorem 1, we have*

$$\text{(56)} \qquad (nb_n)^{-1/2}[S_n(t) - \mathbb{E}S_n(t)] \Rightarrow \text{N}(0, \sigma^2(t)\nu_\mathbf{K}).$$

PROOF. First note that

$$\text{(57)} \qquad \sup_{i\in\mathcal{N}_n(t)} \|\mathcal{P}_{i-k}(e_i - \check{e}_i)\| \leq \sup_{i\in\mathcal{N}_n(t)} (\|\mathcal{P}_{i-k}e_i\| + \|\mathcal{P}_{i-k}\check{e}_i\|) \leq 2\delta_F(k-1).$$

Let $n$ be sufficiently large and $\epsilon > 0$ be sufficiently small. Then for $i \in \mathcal{N}_n(t)$,

$$\|\mathcal{P}_{i-k}(e_i - \check{e}_i)\|$$
$$\leq \|e_i - \check{e}_i\|$$
$$\leq \|(e_i - \check{e}_i) I\{|\zeta_i(i/n) - \zeta(t)| \geq \epsilon\}\|$$

$$+ \|(e_i - \check{e}_i)I\{|\zeta_i(i/n) - \zeta_i(t)| < \epsilon\}\|$$

$$\leq \|I\{|\zeta_i(i/n) - \zeta_i(t)| \geq \epsilon\}\| + \|I\{|\zeta_i(t) - Q(t)| < Q'(t)b_n + \epsilon\}\|$$

$$\leq \epsilon^{-q/2}\|\zeta_i(i/n) - \zeta(t)\|_q^{q/2} + (2Q'(t)b_n + 2\epsilon)^{1/2} f^*(t,\epsilon),$$

where $f^*(t,\epsilon) = \sup\{f(t,x) : x \in (Q(t) - Q'(t)b_n - \epsilon, Q(t) + Q'(t)b_n + \epsilon)\}$. Let $\epsilon = b_n^{2/3}$. By (A2), (A3) and the above inequality, we have

(58) $\quad \|\mathcal{P}_{i-k}(e_i - \check{e}_i)\| \leq C[|i/n - t|^{q/2}/b_n^{q/3} + (b_n^{2/3})^{1/2}] \leq C b_n^{\min(1/3, q/6)}$

for all $i \in \mathcal{N}_n(t)$. Hence by (57) and (58), we obtain

$$\|\mathcal{P}_{i-k}(e_i - \check{e}_i)\| \leq \min\{2\delta_F(k-1), C b_n^{\min(1/3, q/6)}\}.$$

Note that $b_n \to 0$. Using the same argument for (42), we have

(59) $\qquad \|S_n(t) - \mathbb{E}S_n(t) - \check{S}_n(t)\| = o((nb_n)^{1/2}).$

Now by Lemma 4, (56) holds. □

PROOF OF THEOREM 1. Clearly Proposition 5 implies (6). By Proposition 5, we have

(60) $\quad f(t, Q(t))\mu_{\mathbf{K}}(nb_n)^{1/2}\hat{\theta}(t) - \dfrac{[S_n(t) - \mathbb{E}S_n(t)]}{(nb_n)^{1/2}} - \dfrac{\mathbb{E}S_n(t)}{(nb_n)^{1/2}} = o_{\mathbb{P}}(1).$

By (A1), (A2) and Taylor's expansions, since $\alpha = F(i/n, Q(i/n))$, we have

(61)
$$\mathbb{E}S_n(t) = \sum_{i=1}^{n}[F(i/n, Q(i/n)) - F(i/n, Q(t) + Q'(t)(i/n - t))]\mathbf{v}_{i,n}(t)$$
$$= nb_n^3 f(t, Q(t))Q''(t)(1 + o(1))(\mu_2, 0)^\top/2 + o(nb_n^3).$$

By (60), (61) and Proposition 6, we have (8) in view of $nb_n^5 = O(1)$ and

(62) $\quad (nb_n)^{1/2} f(t, Q(t))[\mu_{\mathbf{K}}\hat{\theta}(t) - b_n^2 Q''(t)(\mu_2, 0)^\top/2] \Rightarrow N(0, \sigma^2(t)\nu_{\mathbf{K}}).$ □

PROOF OF COROLLARY 1. By (6), we have

$$\sqrt{nb_n}[\hat{Q}_\alpha(t) - Q_\alpha(t)] - \dfrac{T_{\alpha,n}(t)}{\sqrt{nb_n}f(t, Q_\alpha(t))}[1 + O((nb_n)^{-1})] = o_{\mathbb{P}}(1).$$

Applying the above relation with $\alpha = 0.75$ and $\alpha = 0.25$, we have (9) by using the arguments in the proofs of Theorem 1 and Proposition 6. □

PROOF OF PROPOSITION 1. We shall first prove (ii). A careful check of the proofs of Lemmas 1, 2 and Proposition 5 implies that they are also valid for $t = 0$, and (48) becomes

(63) $\qquad f(0, Q(0))\Sigma_n(0)\hat{\theta}(0) - S_n(0) = o_{\mathbb{P}}((nb_n)^{1/2}).$

NONSTATIONARY QUANTILE 25Also Lemmas 3, 4 and Proposition 6 are still valid with $\nu_{\mathbf{K}}$ therein replaced by the $2 \times 2$ matrix $\nu_{\mathbf{K}}^0 = (\varphi_{|j-j'|})_{1 \le j, j' \le 2}$, where $\varphi_j = \int_0^1 u^j K^2(u)\,du$. By (A1), (A2) and Taylor's expansions, since $K \in \mathcal{K}$, with elementary calculations, (61) becomes

$$\mathbb{E}S_n(0) = \sum_{i=1}^{\lfloor nb_n \rfloor} [F(i/n, Q(i/n)) - F(i/n, Q(0) + Q'(0)i/n)]\mathbf{v}_{i,n}(0)$$

$$= \sum_{i=1}^{n} [Q(i/n) - (Q(0) + Q'(0+)i/n)]f(i/n, Q(i/n))$$

(64)
$$\times (1+o(1))K_{b_n}(i/n)\mathbf{v}_{i,n}(0)$$

$$= \frac{nb_n^3}{2} f(0, Q(0)) Q''(0+)(1 + o(1)) \begin{pmatrix} \int_0^1 u^2 K(u)\,du \\ \int_0^1 u^3 K(u)\,du \end{pmatrix}$$

Additionally, since $K \in \mathcal{K}$, $\Sigma_n(0)/(nb_n) \to \Gamma$, where $\Gamma$ is a $2 \times 2$ matrix with its $(j,j')$th entry being $\int_0^1 u^{|j-j'|} K(u)\,du$, $1 \le j, j' \le 2$. Hence, by (63), $f(0, Q(0))\hat{\theta}(0) - \Sigma_n^{-1}(0)S_n(0) = o_{\mathbb{P}}((nb_n)^{-1/2})$. As in the proof of Theorem 1, the latter relation implies (12) by Slutsky's theorem and (64), since $(nb_n)^{-1/2}[S_n(0) - \mathbb{E}S_n(0)] \Rightarrow \mathrm{N}(0, \sigma^2(0)\nu_{\mathbf{K}}^0)$.

We now prove (i). Let $\varphi_n := \sum_{i=1}^n K_{b_n}(i/n) = (nb_n)/2 + O(1)$ and $\Delta_n = \sum_{i=1}^n \psi(X_i - Q(0))K_{b_n}(i/n)$. Using the argument in the proof of Proposition 5, we have the following analogous version of (63):

(65) $$f(0, Q(0))\varphi_n(\bar{Q}(0) - Q(0)) - \Delta_n = o_{\mathbb{P}}((nb_n)^{1/2}).$$

Similarly as Proposition 6, the CLT $[\Delta_n - \mathbb{E}(\Delta_n)]/\sqrt{nb_n} \Rightarrow \mathrm{N}(0, \sigma^2(0)\varphi_0)$ holds. As (64),

$$\mathbb{E}(\Delta_n) = \sum_{i=1}^{\lfloor nb_n \rfloor} [F(i/n, Q(i/n)) - F(i/n, Q(0))]K_{b_n}(i/n)$$

$$= \sum_{i=1}^{n} [Q(i/n) - Q(0)]f(i/n, Q(i/n))(1+o(1))K_{b_n}(i/n)$$

$$= nb_n^2 f(0, Q(0)) Q'(0+)(1+o(1))\mu_1.$$

Hence (11) holds by elementary manipulations. $\square$

We will only prove Theorem 3 in this section since Theorem 2 can be proved by using the same technique. Let

$$M_n(t, \theta) = \sum_{i=1}^n \{\psi(Y_i(t) - \theta^\top \mathbf{z}_{i,n}(t))$$

(66)
$$- \mathbb{E}[\psi(Y_i(t) - \theta^\top \mathbf{z}_{i,n}(t))|\Upsilon_{i-1}]\}\mathbf{v}_i(t),$$

where $\mathbf{v}_i(t) = (v_{i,1}(t), v_{i,2}(t))^\top = K_{b_n}(i/n - t)(1, (i/n - t)/b_n)^\top$, and

(67)
$$N_n(t, \theta) = \sum_{i=1}^{n} \{\mathbb{E}[\psi(Y_i(t) - \theta^\top \mathbf{z}_{i,n}(t))|\Upsilon_{i-1}]$$
$$- \mathbb{E}\psi(Y_i(t) - \theta^\top \mathbf{z}_{i,n}(t))\}\mathbf{v}_i(t).$$

Then
$$S_n(t, \theta) - \mathbb{E}S_n(t, \theta) = M_n(t, \theta) + N_n(t, \theta).$$

Note that the summands of $M_n(t, \theta)$ form a martingale difference sequence while summands of $N_n(t, \theta)$ are differentiable with respect to $t$ and $\theta$; see Lemmas 5 and 6 below. Decomposing $S_n(t, \theta) - \mathbb{E}S_n(t, \theta)$ into a martingale part and a differentiable part is useful for studying its oscillatory behavior. The technique in Wu (2007a) is useful.

PROOF OF THEOREM 3. Let $\mathcal{S}_n(t, \theta) = \mathbb{E}S_n(t, \theta)$ and $\mathcal{S}_n(t) = \mathbb{E}S_n(t)$. For $\gamma_n \to 0$, we have

$$\sup_{t \in T_n, |\theta| \leq \gamma_n} |S_n(t, \theta) - \mathcal{S}_n(t, \theta) - [S_n(t) - \mathcal{S}_n(t)]|$$
$$= O_\mathbb{P}((nb_n)^{1/2}(\gamma_n^{1/2} \log n + \gamma_n b_n^{-1/2}) + n^{-3})$$

in view of Lemmas 5 and 6. By Lemma 7 below, we have

(68)
$$\sup_{t \in T_n} |S_n(t, \hat{\theta}(t)) - \mathcal{S}_n(t, \hat{\theta}(t)) - [S_n(t) - \mathcal{S}_n(t)]|$$
$$= O_\mathbb{P}((nb_n)^{1/2}(\pi_n^{1/2} \log n + \pi_n b_n^{-1/2})).$$

Elementary calculations using (B1), (B3) and Lemma 7 show that

(69) $\sup_{t \in T_n} |\mathcal{S}_n(t, \hat{\theta}(t)) - \mathcal{S}_n(t) + nb_n f(t, Q(t))\mu_\mathbf{K} \hat{\theta}(t)| = O_\mathbb{P}(nb_n^2 \pi_n + nb_n \pi_n^2).$

Since $(nb_n)^{1/2}\pi_n b_n^{-1/2} = O(nb_n \pi_n^2)$, by (68), (69) and Lemma 8, Theorem 3 follows. □

LEMMA 5. Let $(\gamma_n)_{n \in \mathbb{N}}$ be positive with $\gamma_n \to 0$. Under conditions of Theorem 3, we have

(70) $\sup_{t \in [0,1], |\theta| \leq \gamma_n} |M_n(t, \theta) - M_n(t, 0)| = O_\mathbb{P}((nb_n \gamma_n)^{1/2} \log n + n^{-3}).$

PROOF. Let $\mathcal{N}_n(t-) = (0, nt] \cap \mathcal{N}_n(t)$ and $\mathcal{N}_n(t+) = (nt, n] \cap \mathcal{N}_n(t)$. Define

$$\eta_i(t, \theta) = \psi(Y_i(t) - \theta^\top \mathbf{z}_{i,n}(t)) - \mathbb{E}[\psi(Y_i(t) - \theta^\top \mathbf{z}_{i,n}(t)) | \Upsilon_{i-1}].$$

Write $M_{n,1}(t, \theta) = \sum_{i=1}^n \eta_i(t, \theta) v_{i,1}(t)$, $M_{n,2}(t, \theta) = \sum_{i \in \mathcal{N}_n(t+)} \eta_i(t, \theta) v_{i,2}(t)$ and $M_{n,3}(t, \theta) = \sum_{i \in \mathcal{N}_n(t-)} \eta_{i,\theta}(t) v_{i,2}(t)$. Since $M_n(t, \theta) = (M_{n,1}(t, \theta), \sum_{j=2}^3 M_{n,j}(t, \theta))^\top$, it suffices to prove (70) with $M_n$ therein replaced by $M_{n,j}$ for $j = 1, 2, 3$. For presentation clarity we consider $j = 2$. The other two cases can be similarly dealt with.

For a real sequence $(g_n)_{n \in \mathbb{N}} \to \infty$ with $g_n \geq 3$ for all $n$, define $u_n = (nb_n\gamma_n)^{1/2} g_n / \log g_n$, $\phi_n = (nb_n\gamma_n)^{1/2} g_n \log n$, $a_i(t, \theta) = \psi(Y_i(t) - \theta^\top \mathbf{z}_{i,n}(t)) \times v_{i,2}(t)$, and

$$A_n(t) = \max_{i \in \mathcal{N}_n(t+)} \sup_{|\theta| \leq \gamma_n} |a_i(t, \theta) - a_i(t, 0)|,$$

$$U_n(t) = \sum_{i \in \mathcal{N}_n(t+)} \mathbb{E}\{[\psi(Y_i(t) + |\mathbf{z}_{i,n}(t)|\gamma_n) - \psi(Y_i(t) - |\mathbf{z}_{i,n}(t)|\gamma_n)]^2 | \Upsilon_{i-1}\} v_{i,2}^2(t).$$

By the monotonicity of $\psi(\cdot)$,

(71) $$\sup_{|\theta| \leq \gamma_n} \sum_{i \in \mathcal{N}_n(t+)} \mathbb{E}[(\eta_i(t, \theta) - \eta_i(t, 0))^2 v_{i,2}^2 | \Upsilon_{i-1}] \leq U_n(t).$$

By (B2), (B3) and since $\gamma_n \to 0$, it is easy to see that, for sufficiently large $n$,

$$\mathbb{E}\left[\sup_{t \in [0,1]} U_n(t)\right] \leq 3nb_n\gamma_n c_0 \mathbb{E}\left[\sup_{(t,x) \in [0,1] \times \mathbb{R}} F_1(t, x, \Upsilon_{i-1})\right] \leq Cnb_n\gamma_n,$$

where $c_0 = \int_{-1}^1 K^2(u) u^2 (1 + u^2)^{1/2} du < \infty$. Hence by Markov's Inequality,

(72) $$\mathbb{P}\left[\sup_{t \in [0,1]} U_n(t) \geq u_n^2\right] \leq u_n^{-2} \mathbb{E}\left[\sup_{t \in [0,1]} U_n(t)\right] = O(g_n^{-1} \log g_n) = o(1).$$

Similarly, using $A_n(t) \leq \sup_{|\theta| \leq \gamma_n} \sum_{i \in \mathcal{N}_n(t+)} |a_i(t, \theta) - a_i(t, 0)|$, we have

(73) $$\mathbb{P}\left[\sup_{t \in [0,1]} A_n(t) \geq u_n\right] = o(1).$$

Let $l = n^9$, $\mathcal{G}_l = \{(k_1/l, k_2/l) : |k_1|, |k_2| \leq n^9; k_i \in \mathbb{Z}, i = 1, 2\} \cap \{[-\gamma_n, \gamma_n] \times [-\gamma_n, \gamma_n]\}$ and $\mathcal{H}_l = \{k_3/l : 0 \leq k_3 \leq n^9; k_3 \in \mathbb{Z}\}$. Then $\mathcal{G}_l \times \mathcal{H}_l$ has at most $27n^{27}$ points.

By the argument in the proof of Lemma 4 in Wu (2007a), we obtain in view of Freedman's (1975) exponential inequality for martingale differences that

$$\mathbb{P}\Big\{\sup_{\theta\in\mathcal{G}_l, t\in\mathcal{H}_l} |M_{n,2}(t,\theta) - M_{n,2}(t,0)| \geq 2\phi_n,$$

(74)
$$\sup_{t\in[0,1]} A_n(t) \leq u_n, \sup_{t\in[0,1]} U_n(t) \leq u_n^2\Big\}$$

$$\leq 27 n^{27} O[\exp\{-\phi_n^2/(4u_n\phi_n + 2u_n^2)\}].$$

Since $4u_n\phi_n \log n = o(\phi_n^2)$ and $2u_n^2 \log n = o(\phi_n^2)$, (72), (73) and (74) imply

(75) $$\lim_{n\to\infty} \mathbb{P}\Big\{\sup_{\theta\in\mathcal{G}_l, t\in\mathcal{H}_l} |M_{n,2}(t,\theta) - M_{n,2}(t,0)| \geq 2\phi_n\Big\} = 0.$$

Next we shall apply a chaining argument. For $\mathbf{x} = (x_1, x_2)^\top \in \mathbb{R}^2$, define $\lfloor\mathbf{x}\rfloor_l := (\lfloor x_1\rfloor_l, \lfloor x_2\rfloor_l)^\top$ and $\lceil\mathbf{x}\rceil_l := (\lceil x_1\rceil_l, \lceil x_2\rceil_l)^\top$, where $\lfloor u\rfloor_l = \lfloor ul\rfloor/l$ and $\lceil u\rceil_l = \lceil ul\rceil/l$. Let $\theta \in [-\gamma_n, \gamma_n] \times [-\gamma_n, \gamma_n]$. Since $\psi(\cdot)$ is nondecreasing, we have

(76) $a_i(t, \lceil\theta\rceil_l) \leq a_i(t, \theta) \leq a_i(t, \lfloor\theta\rfloor_l)$   for all $t \in \mathcal{H}_l$ and $i \in \mathcal{N}_n(t+)$.

Let $L_i = \sup_{(t,x)\in[0,1]\times\mathbb{R}} F_1(t, x, \Upsilon_i)$ and $V_n = \sum_{i=1}^n L_{i-1}$. Let $|s_1|, |s_2| \leq \gamma_n$. By (B2), for all large $n$ and $i \in \mathcal{N}_n(t+)$, $|\mathbb{E}[(\psi(Y_i(t)-s_1) - \psi(Y_i(t)-s_2))|\Upsilon_{i-1}]| \leq L_{i-1}|s_1 - s_2|$. Since $|\theta - \lceil\theta\rceil_l| = O(l^{-1})$, we have

(77) $$\sup_{|\theta|\leq\gamma_n} \sum_{i\in\mathcal{N}_n(t+)} |\mathbb{E}\{[a_i(t,\theta) - a_i(t,\lceil\theta\rceil_l)]|\Upsilon_{i-1}\}| \leq Cl^{-1}V_n.$$

The same inequality holds if we replace $\lceil\theta\rceil_l$ by $\lfloor\theta\rfloor_l$. Therefore, for all $|\theta| \leq \gamma_n$,

$$M_{n,2}(t, \lceil\theta\rceil_l) - M_{n,2}(t, 0) - Cl^{-1}V_n \leq M_{n,2}(t, \theta) - M_{n,2}(t, 0)$$
$$\leq M_{n,2}(t, \lfloor\theta\rfloor_l) - M_{n,2}(t, 0) + Cl^{-1}V_n.$$

Since $\mathbb{E}(V_n) \leq Cn$, $l^{-1}V_n = O_\mathbb{P}(n^{-4})$. Since $g_n \to \infty$ can be arbitrarily slow, we have

(78) $$\sup_{t\in\mathcal{H}_l, |\theta|\leq\gamma_n} |M_{n,2}(t,\theta) - M_{n,2}(t,0)| = O_\mathbb{P}((nb_n\gamma_n)^{1/2}\log n + n^{-4}).$$

Let $t_k = k/l$, $k = 0, \ldots, l$. By the triangle inequality, we have

(79)
$$\sup_{t\in[0,1], |\theta|\leq\gamma_n} |M_{n,2}(t,\theta) - M_{n,2}(t,0)|$$
$$\leq \max_{0\leq k\leq l-1} \sup_{0<t-t_k<l^{-1}, |\theta|\leq\gamma_n} |M_{n,2}(t,\theta) - M_{n,2}(t_k,\theta)|$$
$$+ \sup_{t\in\mathcal{H}_l, |\theta|\leq\gamma_n} |M_{n,2}(t,\theta) - M_{n,2}(t,0)|.$$

By the similar chaining argument as those in the proof of (78) as well as the fact that $F_1(t, x, \Upsilon_i)$ is bounded on $[0,1] \times \mathbb{R}$, we have

$$(80) \quad \max_{0 \leq k \leq l-1} \sup_{0 < t - t_k < l^{-1}, |\theta| \leq \gamma_n} |M_{n,2}(t, \theta) - M_{n,2}(t_k, \theta)| = o_{\mathbb{P}}(n^{-3}).$$

By (80), (78) and (79), the lemma follows. $\square$

LEMMA 6. *Under conditions of Lemma 5, we have*

$$(81) \quad \sup_{t \in [0,1], |\theta| \leq \gamma_n} |N_n(t, \theta) - N_n(t, 0)| = O_{\mathbb{P}}(n^{1/2} \gamma_n).$$

PROOF. Let $\Pi_n = \sup_{t \in [0,1], |\theta| \leq \gamma_n} |N_n(t, \theta) - N_n(t, 0)|$. By (B1), we have

$$N_n(t, \theta) - N_n(t, 0) = \int_0^{\theta_1} \int_0^{\theta_2} \frac{\partial^2 N_n(t, (u,v)^\top)}{\partial u \partial v} \, du \, dv$$

$$+ \int_0^{\theta_1} \frac{\partial N_n(t, (u, 0)^\top)}{\partial u} \, du$$

$$+ \int_0^{\theta_2} \frac{\partial N_n(t, (0, v)^\top)}{\partial v} \, dv$$

$$=: \int_0^{\theta_1} \int_0^{\theta_2} \tilde{N}_1(t, \mathbf{u}) \, du \, dv$$

$$+ \int_0^{\theta_1} \tilde{N}_2(t, u) \, du + \int_0^{\theta_2} \tilde{N}_3(t, v) \, dv,$$

where $\mathbf{u} = (u, v)^\top$. Hence,

$$(82) \quad \Pi_n \leq \int_{-\gamma_n}^{\gamma_n} \int_{-\gamma_n}^{\gamma_n} \sup_{t \in [0,1]} |\tilde{N}_1(t, \mathbf{u})| \, du \, dv + \int_{-\gamma_n}^{\gamma_n} \sup_{t \in [0,1]} |\tilde{N}_2(t, u)| \, du$$

$$+ \int_{-\gamma_n}^{\gamma_n} \sup_{t \in [0,1]} |\tilde{N}_3(t, v)| \, dv.$$

Let $\tau_i = i b_n$ for $i = 0, 1, \ldots, \tilde{b}_n$ and $\tau_i = 1$ for $i = \tilde{b}_n + 1$, where $\tilde{b}_n = \lfloor b_n^{-1} \rfloor$. By the triangle inequality, for $\mathbf{u} \in [-\gamma_n, \gamma_n] \times [-\gamma_n, \gamma_n]$, we have

$$(83) \quad \sup_{t \in [0,1]} |\tilde{N}_1(t, \mathbf{u})| \leq \max_{0 \leq i \leq \tilde{b}_n + 1} |\tilde{N}_1(\tau_i, \mathbf{u})| + \max_{1 \leq i \leq \tilde{b}_n + 1} Z_i(\mathbf{u}), \quad \text{where}$$

$$Z_i(\mathbf{u}) = \sup_{\tau_i - b_n < t < \tau_i} |\tilde{N}_1(t, \mathbf{u}) - \tilde{N}_1(\tau_i, \mathbf{u})|.$$

By (B1) and the differentiability of $Q(t), Q'(t), \mathbf{v}_i(t), \mathbf{z}_{i,n}(t)$, we conclude that $\tilde{N}_1(t, \mathbf{u})$ is differentiable with respect to $t$. Furthermore, simple calculations show that

$$(84) \quad \frac{d\tilde{N}_1(t, \mathbf{u})}{dt} = \sum_{k=2}^{3} \sum_{i \in \mathcal{N}_n(t)} [R_{k,n}(i, t, \mathbf{u}) - \mathbb{E}(R_{k,n}(i, \mathbf{u}, t))] \mathbf{w}_{k,n}(i, t, \mathbf{u}),$$

30     Z. ZHOU AND W. B. WU

where $R_{k,n}(i, \mathbf{u}, t) = F_k(i/n, Q(t) + (i/n - t)Q'(t) + \mathbf{u}^\top \mathbf{z}_{i,n}, \Upsilon_{i-1})$ and $\mathbf{w}_{k,n}(i, t, \mathbf{u})$ satisfies

$$\text{(85)} \quad \sup_{t \in [0,1]} \sup_{|\mathbf{u}| \leq \gamma_n} \sum_{i \in \mathcal{N}_n(t)} |\mathbf{w}_{k,n}(i,t,\mathbf{u})|^2 = O(n/b_n), \quad k = 2, 3.$$

Using (85), (B1) and similar arguments as those in the proof of Lemma 1, we have

$$\text{(86)} \quad \sup_{t \in [0,1]} \sup_{|\mathbf{u}| \leq \gamma_n} \left\| \frac{d\tilde{N}_1(t,\mathbf{u})}{dt} \right\| = O((n/b_n)^{1/2}).$$

By (86), we have for $i = 1, 2, \ldots, \tilde{b}_n + 1$

$$\text{(87)} \quad \|Z_i(\mathbf{u})\| \leq \int_{\tau_i - b_n}^{\tau_i} \left\| \frac{d\tilde{N}_1(s,\mathbf{u})}{dt} \right\| ds = O((nb_n)^{1/2}).$$

Since $\max_{1 \leq i \leq \tilde{b}_n + 1} |Z_i(\mathbf{u})|^2 \leq \sum_{i=1}^{\tilde{b}_n+1} |Z_i(\mathbf{u})|^2$, we have

$$\text{(88)} \quad \left\| \max_{1 \leq i \leq \tilde{b}_n+1} |Z_i(\mathbf{u})| \right\| = O((nb_n)^{1/2}(1/b_n)^{1/2}) = O(n^{1/2}).$$

Similarly, $\|\max_{0 \leq i \leq \tilde{b}_n+1} |\tilde{N}_1(\tau_i, \mathbf{u})|\| = O(n^{1/2})$. By (83), we conclude

$$\left\| \sup_{t \in [0,1]} |\tilde{N}_1(t,\mathbf{u})| \right\| = O((nb_n)^{1/2}(1/b_n)^{1/2}) = O(n^{1/2}).$$

Same inequality holds with $\tilde{N}_1(t,\mathbf{u})$ herein replaced by $\tilde{N}_2(t,u)$ or $\tilde{N}_3(t,v)$. Thus by (82), $\|\Pi_n\| = O((nb_n)^{1/2}\gamma_n(1/b_n)^{1/2}) = O(n^{1/2}\gamma_n)$ and this lemma follows. □

LEMMA 7. *Recall $\hat{\theta}(t) = [\hat{Q}(t) - Q(t), b_n(\hat{Q}'(t) - Q'(t))]^\top$ and $\pi_n = (nb_n)^{-1/2}[\log n + (b_n)^{-1/2} + (nb_n^5)^{1/2}]$. Under conditions of Theorem 3, we have $\sup_{t \in [0,1]} |\hat{\theta}(t)| = O_\mathbb{P}(\pi_n)$.*

PROOF. By (47), it suffices to show that

$$\text{(89)} \quad \sup_{t \in [0,1]} |\hat{\vartheta}(t)| = O_\mathbb{P}((nb_n)^{1/2} \pi_n) := O_\mathbb{P}(\bar{\pi}_n),$$

where we recall (45) for $\hat{\vartheta}(t)$. Define

$$\Lambda_n(t, \vartheta) = \sum_{i=1}^{n} [\rho(Y_i(t) - \vartheta^\top \mathbf{u}_{i,n}(t)) - \rho(Y_i(t))] K_{b_n}(i/n - t),$$

$$\Xi_n(t, \vartheta) = \Lambda_n(t, \vartheta) + \sum_{i=1}^{n} \vartheta^\top \mathbf{u}_{i,n}(t) \psi(Y_i(t)) K_{b_n}(i/n - t)$$

(90)
$$:= \Lambda_n(t, \vartheta) + \vartheta^\top \Psi_n(t),$$

$$\tilde{\Xi}_n(t, \vartheta) = \sum_{i=1}^{n} [\psi(Y_i(t) - \vartheta^\top \mathbf{u}_{i,n}(t)) - \psi(Y_i(t))] \mathbf{u}_{i,n}(t) K_{b_n}(i/n - t),$$

where $\mathbf{u}_{i,n}(t)$ is defined in (46). Let $\gamma_n$ be a positive sequence such that

(91) $$\gamma_n \to \infty \quad \text{and} \quad \gamma_n (nb_n)^{-1/2} \to 0.$$

Using the same martingale part and differentiable part decomposition technique applied in Lemmas 5 and 6, we can show that

(92)
$$\sup_{t \in [0,1]} \sup_{|\vartheta| \leq \gamma_n} |\tilde{\Xi}_n(t, \vartheta) - \mathbb{E}\tilde{\Xi}_n(t, \vartheta)|$$
$$= O_{\mathbb{P}}(\gamma_n^{1/2}(nb_n)^{-1/4} \log n + \gamma_n (nb_n)^{-1/2} b_n^{-1/2})$$

and

(93) $$\sup_{t \in [0,1]} |\Psi_n(t) - \mathbb{E}[\Psi_n(t)]| = O_{\mathbb{P}}(\log n + (b_n)^{-1/2}).$$

On the other hand, simple calculations show that $\sup_{t \in [0,1]} |\mathbb{E}[\Psi_n(t)]| = O_{\mathbb{P}}((nb_n^5)^{1/2})$. Together with (93) we conclude that

(94) $$\sup_{t \in [0,1]} |\Psi_n(t)| = O_{\mathbb{P}}(\bar{\pi}_n).$$

Furthermore, using the similar arguments as those in the proof of Lemma 2 and (49), it is easy to see that for $\gamma_n$ satisfying (91)

(95) $$\sup_{t \in [0,1]} \sup_{|\vartheta| \leq \gamma_n} |\mathbb{E}[\Xi_n(t, \vartheta)] - f(t, Q(t))|\vartheta|^2/2| = o(\gamma_n^2).$$

Note $\Xi_n(t, \vartheta) = -\vartheta^\top \int_0^1 \tilde{\Xi}_n(t, s\vartheta)\, ds$. Relation (92) implies that

$$\sup_{t \in [0,1]} \sup_{|\vartheta| \leq \gamma_n} |\Xi_n(t, \vartheta) - \mathbb{E}\Xi_n(t, \vartheta)| \leq \gamma_n \sup_{t \in [0,1]} \sup_{|\vartheta| \leq \gamma_n} |\tilde{\Xi}_n(t, \vartheta) - \mathbb{E}\tilde{\Xi}_n(t, \vartheta)|$$
$$= o_{\mathbb{P}}(\gamma_n^2),$$

which by (95) implies

(96) $$\sup_{t \in [0,1]} \sup_{|\vartheta| \leq \gamma_n} \left| \Xi_n(t, \vartheta) - \frac{1}{2} f(t, Q(t))|\vartheta|^2 \right| = o_{\mathbb{P}}(\gamma_n^2).$$

For a sequence $\{c_n\}$ with $c_n \to \infty$, let $c_n' = \min\{c_n, (\pi_n)^{-1/2}\}$. Then $c_n' \to \infty$. It is easy to see that $\bar{\gamma}_n = c_n' \bar{\pi}_n$ satisfies (91). On the other hand, if $|\vartheta| = \bar{\gamma}_n$, we have by (94) that $\sup_{t \in [0,1]} |\vartheta^\top \Psi_n(t)| = O_{\mathbb{P}}(c_n' \bar{\pi}_n^2)$. Since $\inf_{t \in [0,1]} f(t, Q_\alpha(t)) > 0$, by (96), we have

(97) $$\mathbb{P}\left[ \inf_{t \in [0,1]} \inf_{|\vartheta| = \bar{\gamma}_n} \Lambda_n(t, \vartheta) \leq 0 \right] \to 0.$$

Note that $\Lambda_n(t, \vartheta)$ is convex in $\vartheta$ and $\Lambda_n(t, 0) = 0$. Hence for any $\vartheta$ such that $|\vartheta| > \bar{\gamma}_n$, we have $\Lambda_n(t, \vartheta) \geq (|\vartheta|/\bar{\gamma}_n)\Lambda_n(t, \bar{\gamma}_n\vartheta/|\vartheta|)$. By (97), $\mathbb{P}[\inf_{t \in [0,1]} \inf_{|\vartheta| \geq \bar{\gamma}_n} \Lambda_n(t, \vartheta) \leq 0] \to 0$. By the definition of $\hat{\vartheta}(t)$, we have $\mathbb{P}[\sup_{t \in [0,1]} |\hat{\vartheta}(t)| \geq c'_n \bar{\pi}_n] \to 0$. Since $c_n \to \infty$ can be arbitrarily slow, (89) follows. $\square$

REMARK 1. Note that under conditions of Theorem 3, $\pi_n \to 0$. Hence Lemma 7 implies that $\hat{Q}(t)$ is a uniformly consistent estimator of $Q(t)$ on $[0, 1]$. Additionally, if $nb_n^4 \to \infty$, then $\hat{Q}'(t)$ is also uniformly consistent for $Q'(t)$ on $[0, 1]$.

LEMMA 8. *Under* (B2), $\sup_{t \in [0,1]} |S_n(t, \hat{\theta}(t))| = O_{\mathbb{P}}(1)$.

PROOF. From the proof of Corollary 2 in Wu (2007a) or Lemma A.2 in Ruppert and Carroll (1980), it is easy to see that

$$(98) \quad \sup_{t \in [0,1]} |S_n(t, \hat{\theta}(t))| \leq \sup_{t \in [0,1]} \sum_{i=1}^{n} |\mathbf{v}_i(t)| I\{X_i = \hat{\mathbf{Q}}(t)^\top \mathbf{z}_{i,n}(t)\},$$

where $\hat{\mathbf{Q}}(t) = (\hat{Q}(t), b_n \hat{Q}'(t))$. Note that

$$\sup_{t \in [0,1]} \sum_{i=1}^{n} |\mathbf{v}_i(t)| I\{X_i = \hat{\mathbf{Q}}(t)^\top \mathbf{z}_{i,n}(t)\}$$

$$\leq C \sup_{t \in [0,1]} \sum_{i=1}^{n} I\{X_i = \hat{\mathbf{Q}}(t)^\top \mathbf{z}_{i,n}(t)\}$$

$$= C \sup_{t \in [0,1]} \sum_{i=1}^{n} I\{X_i = \hat{Q}(t) - t\hat{Q}'(t) + \hat{Q}'(t)i/n\}$$

$$\leq C \sup_{\mathbf{a} \in \mathbb{R}^2} \sum_{i=1}^{n} I\{X_i = a_1 + a_2 i/n\},$$

where $\mathbf{a} = (a_1, a_2)^\top$. Following Babu (1989) and using (B2), we see that

$$\sup_{\mathbf{a} \in \mathbb{R}^2} \sum_{i=1}^{n} I\{X_i = a_1 + a_2 i/n\} = O_{\mathbb{P}}(1).$$

Therefore this lemma follows. $\square$

PROOF OF THEOREM 4. It suffices to prove (20). By (17), we have

$$\sqrt{nb_n}\left[\check{Q}_\alpha(t) - \sum_{i=1}^{n} w_n(t,i)Q_\alpha(i/n)\right] - \sum_{i=1}^{n} \frac{w_n(t,i)T_{\alpha,n}(i/n)}{f(i/n, Q_\alpha(i/n))\sqrt{nb_n}} = o_{\mathbb{P}}(1).$$

Since $\sup_{|s-t|\leq \bar{b}_n} |T_{\alpha,n}(s)| = O_{\mathbb{P}}(nb_n\varrho_n)$, $\sum_{i=1}^n w_n(t,i)Q_\alpha(i/n) - Q_\alpha(t) = O(\bar{b}_n^2)$ and, by (B1)–(B3), $\sup_{|i/n-t|\leq \bar{b}_n} |1/f(i/n, Q_\alpha(i/n)) - 1/f(t, Q_\alpha(t))| = O(\bar{b}_n)$, we have

$$\sqrt{nb_n} f(t, Q_\alpha(t))[\check{Q}_\alpha(t) - Q_\alpha(t)] - \sum_{i=1}^n w_n(t,i) T_{\alpha,n}(i/n)/\sqrt{nb_n} = o_{\mathbb{P}}(1).$$

(99)

Basic manipulations similar as those in the proof of Proposition 6 show that

$$\left\| \sum_{i=1}^n w_n(t,i) \frac{T_{\alpha,n}(i/n) - T_{\alpha,n}(t) - \mathbb{E}\{T_{\alpha,n}(i/n) - T_{\alpha,n}(t)\}}{\sqrt{nb_n}} \right\| = o(\bar{b}_n/b_n).$$

(100)

Since $K \in \mathcal{K}$, it is easily seen that $\sup_{|i/n-t|\leq \bar{b}_n} |\mathbb{E}\{T_{\alpha,n}(i/n) - T_{\alpha,n}(t)\}| = o(\sqrt{nb_n})$. Combining (99) and (100), we have (20). $\square$

OUTLINE OF THE PROOF OF THEOREM 5. From the proofs of Lemmas 5 and 6, we can obtain $\sup_{s\in\mathcal{N}_n(t)} |\hat{Q}_\alpha(s) - Q_\alpha(s)| = O_{\mathbb{P}}(\varrho_n)$. Using the martingale decomposition technique therein, it can be shown that, if we replace the estimates $\hat{Q}(i/n)$ and $\hat{Q}(t)$ in (25) and (26) by the true values $Q(i/n)$ and $Q(t)$, the changes of (25) and (26) will be of order $o_{\mathbb{P}}(1)$. By the property of local stationarity, the same conclusion holds when we further replace $Z_{i,\alpha}$ and $X_i$ in (25) and (26) by $\check{e}_i = \psi_\alpha(\check{X}_i - Q(t))$ and $\check{X}_i = G(t, \Upsilon_i)$, respectively. Applying convergence results for the stationary processes $\{\check{e}_i\}$ and $\{\check{X}_i\}$, Theorem 5 follows. A detailed proof is available upon request. $\square$

**Acknowledgments.** We are grateful to three referees, the associate editor and the editor for their many helpful comments.

Department of Statistics
The University of Chicago
5734 S. University Avenue
Chicago, Illinois 60637
USA
E-mail: zhou@galton.uchicago.edu
wbwu@galton.uchicago.edu